\documentstyle[12pt]{amsart}

\allowdisplaybreaks

\newtheorem{theorem}{Theorem}[section]

\theoremstyle{definition}

\theoremstyle{remark}

\numberwithin{equation}{section}

\renewcommand{\l}{\left}
\newcommand{\r}{\right}
\renewcommand{\a}{\alpha}
\renewcommand{\b}{\beta}
\newcommand{\g}{\gamma}
\newcommand{\f}{\varphi}
\renewcommand{\t}{\theta}
\newcommand{\ep}{\varepsilon}
\newcommand{\om}{\omega}

\begin{document}

\title [Some orthogonal $_8\f_7$-functions]%
{Some orthogonal very-well-poised $_8\f_7$-functions\\
that generalize Askey--Wilson polynomials}

\author{Sergei~K.~Suslov}
\address
{Department of Mathematics, Arizona State University,
Tempe, Arizona 85287--1804}
%{Kurchatov Institute, Moscow 123182, Russia}
%\curraddr{Mathematical Sciences Research Institute,
%1000 Centennial Drive, Berkeley, CA 94720}
\email{suslov{@@}math.la.asu.edu}
%\quad suslov{@@}msri.org}

\thanks{Research at MSRI supported in part by NSF grant \# DMS 90220140.}

\date{July 1, 1997}


\begin{abstract}
In a recent paper Ismail, Masson, and Suslov \cite{Is:Ma:Su2} have
established a continuous orthogonality relation and some other properties
of a $_2\f_1$-Bessel function on a $q$-quadratic grid.
Dick Askey \cite{As} suggested that the ``Bessel-type orthogonality''
found in \cite{Is:Ma:Su2} at the $_2\f_1$-level has really
a general character and can be extended up to the $_8\f_7$-level.
Very-well-poised $_8\f_7$-functions are known as a nonterminating version of
the classical Askey--Wilson polynomials \cite{As:Wi0}, \cite{As:Wi}.
Askey's congecture has been proved in \cite{Su1}.
In the present paper which is an extended version of \cite{Su1}
we discuss in details some properties of the orthogonal $_8\f_7$-functions.
Another type of the orthogonality relation for
a very-well-poised $_8\f_7$-function
was recently found by Askey, Rahman, and Suslov \cite{As:Ra:Su}.
\end{abstract}

\maketitle

\section{Introduction}

The Askey--Wilson polynomials \cite{As:Wi} are
\begin{align}
p_n(x)&= p_n(x; a, b, c, d)  \\
&=a^{-n}\;\l(ab,\, ac,\, ad ;q \r)_n
\ _4\f_3\l(
\begin{array}{cc}
q^{-n}, abcdq^{n-1}, a e^{i\t}, a e^{-i\t} \\
\, \\
ab, ac, ad
\end{array}
;\; q, \, q \r), \notag
\end{align}
where $x=\cos \t$. These polynomials are the most general known
classical orthogonal polynomials (see \cite{An:As}, \cite{As:Wi},
\cite{At:Ra:Su}, \cite{At:Su}, \cite{Ga:Ra}, \cite{Ko:Sw},
and \cite{Ni:Su:Uv}).

The symbol $_4\f_3$ in (1.1) is a special case of
basic hypergeometric series \cite{Ga:Ra},
\begin{equation}
\begin{split}
\ _r\f_s(t):=& \ _r\f_s\l(
\begin{array}{cc}
a_1, a_2, \ldots, a_r \\
\, \\
b_1, b_2, \ldots ,b_s
\end{array}
;\; q, \,t\r) \\
=& \sum^\infty_{n=0}
\frac{(a_1, a_2, \ldots, a_r;q)_n}
{(q, b_1, b_2, \ldots ,b_s;q)_n}
\, \l((-1)^n q^{n(n-1)/2}\r)^{1+s-r}\,t^n.
\end{split}
\end{equation}
The standard notations for the $q$-shifted factorials are
\begin{align}
&(a;q)_n := \prod_{k=1}^n(1 - aq^{k-1}),  \\
&(a_1, a_2, \ldots, a_r;q)_n := \prod_{k=1}^r (a_k;q)_n,
\end{align}
and
\begin{align}
&(a;q)_\infty := \lim_{n\to\infty}(a; q)_n,  \\
&(a_1, a_2, \ldots, a_r;q)_\infty := \prod_{k=1}^r (a_k;q)_\infty
\end{align}
provided $|q|<1$.
For an excellent account on the theory of basic hypergeometric series
see \cite{Ga:Ra}. We shall just mention here that the $_{s+1}\f_s$-series
is called balanced if $qa_1 a_2 \cdots a_{s+1}=b_1 b_2 \cdots b_s$ and
$t=q$; it is a very-well-poised if
$a_2 b_1=a_3 b_2= \ldots =a_{s+1} b_s=q a_1$, and $a_2=q\sqrt{a_1}$,
$a_3=-q\sqrt{a_1}$.

Askey and Wilson found the orthogonality relation
\begin{align}
&\int_{0}^{\pi}\,
\frac{p_n(\cos \t; a, b, c, d)\, p_m(\cos \t; a, b, c, d)\,
\l(e^{2i \t}, e^{-2i \t} ; q\r)_\infty}
{\l(a e^{i \t}, a e^{-i \t}, b e^{i \t}, b e^{-i \t},
c e^{i \t}, c e^{-i \t}, d e^{i \t}, d e^{-i \t}; q \r)_\infty}\, d \t \\
&\quad = \delta_{n m}\; \frac{(q, ab, ac, ad, bc, bd, cd; q)_\infty}
{2 \pi\, (abcd; q)_\infty}  \notag \\
&\qquad \times \frac{\l(1-abcd q^{-1}\r)\,
\l(q, ab, ac, ad, bc, bd, cd; q\r)_n}
{\l(1-abcd q^{2n-1}\r)\, \l(abcd q^{-1}; q\r)_n} \notag.
\end{align}
In the fundamental memoir \cite{As:Wi}, they studied in details
many other properties of these polynomials.

As is well known, the Askey--Wilson polynomials and their special and
limiting cases are the simplest and the most important orthogonal solutions
of a difference equation of hypergeometric type on nonuniform lattices
(see, for example, \cite{As:Wi}, \cite{At:Ra:Su}, \cite{At:Su},
\cite{At:Su1}, \cite{Ga:Ra}, \cite{Ni:Su:Uv}, and \cite{Su}).
Recently Ismail, Masson, and Suslov have found another type of
orthogonal solutions of this difference equation
\cite{Is:Ma:Su2}--\cite{Is:Ma:Su3}.
They considered the $_2\f_1$-function,
\begin{align}
&J_\nu(z, r)=
\widetilde{J}_\nu(x(z), r|q) \\
&\quad:= \l(\frac{r}{2}\r)^\nu\;
\frac{\l( q^{\nu+1}, \, -r^2/4 ; q \r)_\infty}{(q; q)_\infty} \;
\ _2\f_1 \l(
\begin{array}{cc}
q^{(\nu+1)/2} e^{i\t}, \;
q^{(\nu+1)/2} e^{-i\t} \\
\, \\
q^{\nu+1}
\end{array}
;\; q, \, - \frac{r^2}{4} \r), \notag
\end{align}
as a $q$-analog on a $q$-quadratic grid of the Bessel function \cite{Wa},
\begin{equation}
J_{\nu}(x)= \l(\frac{x}{2}\r)^\nu\;
\sum_{n=0}^\infty \; \frac{\l(-x^2/4\r)^n}{n!\, \Gamma(\nu + n +1)}.
\end{equation}
Ismail, Masson, and Suslov established
the following orthogonality property for the $q$-Bessel function,
\begin{align}
&\int_0^\pi
\widetilde{J}_\nu(\cos \t, r)\, \widetilde{J}_\nu(\cos \t, r')\,
\frac{\l(e^{2i\t}, \, e^{-2i\t}; q\r)_\infty}
{\l(q^\a e^{i\t}, q^\a e^{-i\t},
q^{1-\a}e^{i\t}, q^{1-\a}e^{-i\t}; q\r)_\infty} \\
&\qquad\times \l(q^{(\nu+1)/2}e^{i\t}, q^{(\nu+1)/2}e^{-i\t},
q^{(\nu+1)/2}e^{i\t}, q^{(\nu+1)/2}e^{-i\t} ;q\r)_\infty^{-1}\,d\t= 0 \notag
\end{align}
if $ r \ne r'$,
\begin{align}
&\int_0^\pi
\l(\widetilde{J}_\nu(\cos \t, r)\r)^2\,
\frac{\l(e^{2i\t}, \, e^{-2i\t}; q\r)_\infty}
{\l(q^\a e^{i\t}, q^\a e^{-i\t},
q^{1-\a}e^{i\t}, q^{1-\a}e^{-i\t}; q\r)_\infty} \\
&\qquad\times \l(q^{(\nu+1)/2}e^{i\t}, q^{(\nu+1)/2}e^{-i\t},
q^{(\nu+1)/2}e^{i\t}, q^{(\nu+1)/2}e^{-i\t} ;q\r)_\infty^{-1}\,d\t \notag\\
&\qquad\qquad=
\dfrac{-4\pi(1-q)q^{-(\nu+1)/2}}
{\l(q, q^{(\nu+1)/2+\a}, q^{(\nu+1)/2-\a+1}; q\r)_\infty^2}\,
\dfrac{\partial \widetilde{J}_\nu(x(\a),r)}{\partial r^2}\;
\dfrac{\nabla \widetilde{J}_\nu(x(\a),r)}{\nabla x(\a)} \notag
\end{align}
if $r = r'$.
Here $r$ and $r'$ are two roots of the equation
\begin{equation}
J_\nu(\a, r)=J_\nu(\a, r')=0,
\end{equation}
and $J_\nu(z, r)=\widetilde{J}_\nu(x(z),r)$,
$x(z)=\dfrac{1}{2}\l(q^z+q^{-z}\r)$, ($x=\cos \t$, if $q^z=e^{i\t}$),
$\text{Re }\nu >-1$, and $0<\text{Re }\a<1$.
(See \cite{Is:Ma:Su2}--\cite{Is:Ma:Su3} for more details.)
This is a $q$-version of the orthogonality relation
for the classical Bessel function
\begin{equation}
\int_0^1 x\, J_\nu(r x)\, J_\nu(r'x) \, dx =
\begin{cases}
0& \text{if $ r \ne r'$}, \\
\dfrac{1}{2}\, \l(J_{\nu+1}(r)\r)^2 & \text{if $r = r'$},
\end{cases}
\end{equation}
under the conditions $J_\nu(r)=J_\nu(r')=0$ \cite{Wa}.

Another example of an orthogonality relation of
this type has been discovered recently. The following functions
\begin{align}
C(x):=& C_q(x; \om) \\
= &\frac{(-\om^2; q^2)_\infty}{(-q\om^2; q^2)_\infty}\;
\ _2\f_1\l(
\begin{array}{cc}
-qe^{2i\t}, \, -qe^{-2i\t}\\ \, \\ q
\end{array}
;\; q^2, \, -\om^2 \r) \notag
\end{align}
and
\begin{align}
S(x):=& S_q(x; \om) \\
= &\frac{(-\om^2; q^2)_\infty}{(-q\om^2; q^2)_\infty}\;
\frac{2q^{1/4}}{1-q} \; \om \notag \\
&\times \cos \t \ _2\f_1\l(
\begin{array}{cc}
-q^2e^{2i\t}, \, -q^2e^{-2i\t}\\ \, \\ q^3
\end{array}
;\; q^2, \, -\om^2 \r) \notag
\end{align}
were discussed in
\cite{At:Su} and \cite{Is:Zh} as analogs of
$\cos \om x$ and $\sin \om x$ on a $q$-quadratic lattice, respectively.
Ismail and Zhang \cite{Is:Zh} expanded the corresponding basic exponential
function in terms of ``$q$-spherical harmonics'' and later, together with
Rahman \cite{Is:Ra:Zh}, they extended this $q$-analog of the expansion
formula of the plane wave from $q$-ultraspherical polynomials to
continuous $q$-Jacobi polynomials.
``Addition" theorems for the basic trigonometric functions (1.14)--(1.15)
where found in \cite{Su0}.
Bustoz and Suslov \cite{Bu:Su} have established
an orthogonality property,
\begin{equation}
\int_0^\pi
C(\cos \t;\om)\, C(\cos \t;\om')\,
\frac{\l(e^{2i\t}, \, e^{-2i\t}; q\r)_\infty}
{\l(q^{1/2} e^{2i\t}, q^{1/2} e^{-2i\t}; q \r)_\infty}
\,d\t= 0,
\end{equation}
\begin{equation}
\int_0^\pi
S(\cos \t;\om)\, S(\cos \t;\om')\,
\frac{\l(e^{2i\t}, \, e^{-2i\t}; q\r)_\infty}
{\l(q^{1/2} e^{2i\t}, q^{1/2} e^{-2i\t}; q \r)_\infty}
\,d\t= 0,
\end{equation}
\begin{equation}
\int_0^\pi
C(\cos \t;\om)\, S(\cos \t;\om')\,
\frac{\l(e^{2i\t}, \, e^{-2i\t}; q\r)_\infty}
{\l(q^{1/2} e^{2i\t}, q^{1/2} e^{-2i\t}; q \r)_\infty}
\,d\t= 0,
\end{equation}
and
\begin{align}
&\int_0^\pi
C^2(\cos \t;\om)\,
\frac{\l(e^{2i\t}, \, e^{-2i\t}; q \r)_\infty}
{\l(q^{1/2} e^{2i\t}, q^{1/2} e^{-2i\t}; q \r)_\infty}\, d\t \\
&\quad=\int_0^\pi
S^2(\cos \t;\om)\,
\frac{\l(e^{2i\t}, \, e^{-2i\t}; q \r)_\infty}
{\l(q^{1/2} e^{2i\t}, q^{1/2} e^{-2i\t}; q \r)_\infty} \, d\t \notag\\
&\quad=\pi \frac{\l(q^{1/2}, -q^{1/2} \om^2; q \r)_\infty}
{\l(q, -\om^2; q \r)_\infty}\,
\frac{\l(-\om^2; q^2 \r)_\infty}
{\l(-q\om^2; q^2 \r)_\infty} \notag\\
&\qquad\times
\ _2\f_1\l(
\begin{array}{cc}
-q^{1/2}, \, -\om^2 \\ \, \\ -q^{1/2} \om^2
\end{array}
;\; q, \, q \r). \notag
\end{align}
Here $\om$ and $\om'$ are different solutions of the equation
\begin{equation}
S\l(\frac{1}{2}\l(q^{1/4}+q^{-1/4}\r); \om \r)=0.
\end{equation}
Bustoz and Suslov introduced the corresponding $q$-Fourier series
and established several important facts about the basic trigonometric
system and the $q$-Fourier series.

Dick Askey \cite{As} has suggested that the orthogonality
relation (1.10)--(1.11)
can be extended to the level of very-well-poised $_8\f_7$-functions.
The author was able to prove his conjecture in \cite{Su1}.
Our main objective in this paper is to study the new orthogonality property
for $_8\f_7$-function in details.
This article is an extended version of the short paper \cite{Su1},
originally submitted as a Letter,
we include proofs of the main results etablished in \cite{Su1}
to make this work as self-contained as possible.

The paper is organized as follows.
In Section 2 we consider difference equation of hypergeometric type
on a $q$-quadratic grid and discuss a solution of this equation
in terms of a very-well-poised $_8\f_7$-function.
In the next two sections we derive a continuous orthogonality property
for this $_8\f_7$-function.
Section 5 is devoited to the investigation of zeros and asymptotics
of this function and in Section 6 we evaluate the normalization
constants in the orthogonality relation for the $_8\f_7$.
In Section 7 we find an analog of the Wronskian determinant
for solutions of difference equation of hypergeometric type.
Some special and limiting cases of our orthogonality relation
are discussed in Section 8.
We close this paper by estimating the number of zeros on the basis
of Jensen's theorem in Section 9.

\section{Difference Equation and Its $_8\f_7$-Solutions}

Let us consider a {\it difference equation of hypergeometric type\/}
\begin{equation}
\sigma(z) \,\frac{\Delta}{\nabla x_1(z)}\l(\frac{\nabla u(z)}{\nabla x(z)}\r)
+\tau (z) \, \frac{\Delta u(z)}{\Delta x(z)}
+ \lambda \, u(z) = 0,
\end{equation}
on a $q$-quadratic lattice $x(z) = \dfrac{1}{2}\l(q^z + q^{-z}\r)$
with $x_1(z) = x\l(z + \frac{1}{2}\r)$ and
$\Delta f(z) = \nabla f(z+1) = f(z+1)-f(z)$.
Here, in the most general case,
\begin{align}
\sigma(z) &= q^{-2z}\,\l(q^z-a\r) \l(q^z-b\r) \l(q^z-c\r) \l(q^z-d\r), \\
\tau(z) &= \frac{\sigma(-z)-\sigma(z)}{\nabla x_1(z)} \notag \\
&= \frac{2q^{1/2}}{1-q}
\l(abc+abd+acd+bcd-a-b-c-d+2(1-abcd)\,x \r), \\
\lambda &= \lambda_{\nu}=\frac{4q^{3/2}}{(1-q)^2}\;\l(1-q^{-\nu}\r)
\l(1-abcd\, q^{\nu-1} \r).
\end{align}
Equation (2.1) can also be rewritten in self-adjoint form,
\begin{equation}
\frac{\Delta}{\nabla x_1(z)}
\l(\sigma(z)\,\rho(z)\, \frac{\nabla u(z)}{\nabla x(z)} \r)
+ \lambda\, \rho(z)\, u(z)=0,
\end{equation}
where $\rho(z)$ is a solution of the {\it Pearson equation\/},
\begin{equation}
\Delta \l(\sigma(z)\, \rho(z) \r) =\tau(z)\, \nabla x_1(z).
\end{equation}
See \cite{Ni:Su:Uv} and \cite{Su} for details.

Methods of solving of the difference equation of hypergeometric type (2.1)
were discussed in \cite{At:Su}, \cite{Gru:Hai}, and \cite{Su}.
As is well known, there are different kinds of solutions of this equation.
For integer values of the parameter $\nu=n=0, 1, 2, \ldots ,$
the famous solutions of (2.1) are the Askey--Wilson
$_4\f_3$-polynomials (see, for example, \cite{As:Wi0}, \cite{As:Wi},
\cite{At:Su1}, and \cite{Ga:Ra}).
For arbitrary values
of this parameter solutions of (2.1) can be written in terms of
$_8\f_7$-functions \cite{At:Su}, \cite{Is:Ra}, and \cite{Su}.
Let us choose the following solution,
$u_\nu(z)=u_\nu\l(x(z); a, b, c; d\r)$,
such that
\begin{align}
&u_{\nu}(x; a, b, c; d) \notag \\
&=\frac
{\l(qa/d, bcq^{\nu}, q^{1-\nu+z}/d, q^{1-\nu-z}/d; q\r)_\infty}
{\l(q^{1-\nu}a/d, bc, q^{1+z}/d, q^{1-z}/d; q\r)_\infty}\\
&\quad\times \ _8\f_7\l(
\begin{array}{cc}
\dfrac{aq^{-\nu}}{d},
q \sqrt{\dfrac{aq^{-\nu}}{d}},-q \sqrt{\dfrac{aq^{-\nu}}{d}},
q^{-\nu}, \dfrac{q^{1-\nu}}{bd}, \dfrac{q^{1-\nu}}{cd}, a q^z, aq^{-z} \\
\, \\
\sqrt{\dfrac{aq^{-\nu}}{d}}, -\sqrt{\dfrac{aq^{-\nu}}{d}},
\dfrac{aq}{d}, ab, ac, \dfrac{q^{1-\nu+z}}{d}, \dfrac{q^{1-\nu-z}}{d}
\end{array}
;\; q, \, bcq^{\nu} \r) \notag\\
&=\frac{\l(bcq^{\nu}, q^{1-\nu}/ad;q \r)_\infty}
{\l(bc, q/ad;q \r)_\infty} \;
\ _4\f_3\l(
\begin{array}{cc}
q^{-\nu}, abcdq^{\nu-1}, a q^z, aq^{-z} \\
\, \\
ab, ac, ad
\end{array}
;\; q, \, q \r) \\
&\quad+
\frac{\l(q^{-\nu}, abcd q^{\nu-1}, qb/d, qc/d;q \r)_\infty}
{\l(ab, ac, bc, ad/q ;q \r)_\infty}  \notag \\
&\quad\times \frac{\l(aq^z, aq^{-z}; q\r)_\infty}
{\l( q^{1+z}/d, q^{1-z}/d; q\r)_\infty}\;
\ _4\f_3\l(
\begin{array}{cc}
q^{1-\nu}/ad, bcq^{\nu}, q^{1+z}/d, q^{1-z}/d \\
\, \\
qb/d, qc/d, q^2/ad
\end{array}
;\; q, \, q \r). \notag
\end{align}
We have used Bailey's formula (III.36) of \cite{Ga:Ra}
to transform (2.7) to (2.8).
Both $_4\f_3$-functions in (2.8) are balanced and converge when
$|q|<1$. The $_8\f_7$-function in (2.7) has a very-well-poised structure,
it converges when $|bcq^{\nu}|<1$; (2.8) provides
an analytic continuation of (2.7).
This $_8\f_7$ can be transformed to another equivalent form by
(III.23)--(III.24) of \cite{Ga:Ra}.

A similar function was earlier discussed by Rahman \cite{Ra} and
in a recent paper \cite{Ra1} he has found a $q$-extension of a product
formula of Watson involving this function.
Rahman has also shown that
\begin{align}
&\lim_{q\to1^-}
u_{\nu}(x; q^{\a+1/2}, q^{\a+3/4}, -q^{\b+1/2}; -q^{\b+3/4}) \\
&\quad=\ _2F_1\l(-\nu, \a+\b+\nu+1; \a+1; \dfrac{1-x}{2}\r). \notag
\end{align}
Thus, the $_8\f_7$-function in (2.7) can be thought as a $q$-extension
of the hypergeometric function of Gauss.

It is easy to see that $u_\nu(x; a, b, c; d)$ defined by (2.7)--(2.8)
is a function in $x=\l(q^z+q^{-z}\r)/2$. Indeed, by (1.3) and (1.4),
\begin{equation*}
\l(\xi q^z, \xi q^{-z}\r)_n
=\prod_{k=0}^{n-1} \l(1-2\xi x q^k + \xi^2 q^{2k}\r),
\end{equation*}
where $\xi =a,\, q/d, \, q^{1-\nu}/d$ and $n=1, 2, 3, \ldots$ or $\infty$.
It is also important to mention that
\begin{align}
&\l(q^{-\nu}, abcdq^{\nu-1}; q\r)_n \notag\\
&=\prod_{k=0}^{n-1}
\l(1-\l(q^{-\nu}+abcdq^{\nu-1}\r)q^k + abcd\; q^{2k-1}\r) \notag\\
&=\prod_{k=0}^{n-1}
\l(1-\l(1+abcdq^{-1}-\frac{(1-q)^2}{4q^{3/2}}\,\lambda_\nu \r)q^k
+abcd\; q^{2k-1}\r) \notag
\end{align}
and
\begin{align}
&\l(bcq^{\nu}, q^{1-\nu}/ad; q\r)_n \notag\\
&=\prod_{k=0}^{n-1}
\l(1-\frac{1}{ad}\,\l(q^{-\nu}+abcdq^{\nu-1}\r) q^{k+1}
+\frac{bc}{ad}\; q^{2k+1}\r) \notag\\
&=\prod_{k=0}^{n-1}
\l(1-\frac{1}{ad}\,\l(1+abcdq^{-1}-\frac{(1-q)^2}{4q^{3/2}}\,\lambda_\nu \r)
q^{k+1}+\frac{bc}{ad}\; q^{2k+1}\r) \notag
\end{align}
by (1.3)--(1.4) and (2.4), where $n=1, 2, 3, \ldots$ or $\infty$.
Therefore, $u_\nu(x; a, b, c; d)$ is, really, a function of $\lambda_\nu$
as well.

The definition of the $u_\nu$ in this paper is the same as in
\cite{Ra} and \cite{Ra1},  but different from one in \cite{Su1}.
This definition emphasizes the symmetry properties.
Function $u_\nu(x; a, b, c; d)$ in (2.7)--(2.8) is obviously
symmetric in $b$ and~$c$. Applying (III.36) and Bailey's transform (III.36)
of \cite{Ga:Ra} once again we obtain
\begin{align}
&u_{\nu}(x; a, b, c; d) \notag \\
&=\frac
{\l(abcq^{\nu+z}, q^{1-\nu+z}/d; q\r)_\infty}
{\l(abcq^z, q^{1+z}/d; q\r)_\infty}\\
&\quad\times
\ _8W_7\l(abcq^{z-1};
q^{-\nu}, abcdq^{\nu-1}, aq^z, bq^z, cq^z
;\; q, \, \dfrac{q^{1-z}}{d} \r) \notag\\
&=\frac{\l(abq^{\nu}, acq^{\nu}, bcq^{\nu},
abcdq^{\nu-1}, q^{1-\nu+z}/d, q^{1-\nu-z}/d;q \r)_\infty}
{\l(ab, ac, bc, abcdq^{2\nu-1}, q^{1+z}/d, q^{1-z}/d;q \r)_\infty} \\
&\quad\times \ _4\f_3\l(
\begin{array}{cc}
q^{-\nu}, q^{1-\nu}/{ad}, q^{1-\nu}/{bd}, q^{1-\nu}/{cd} \\
\, \\
q^{2-2\nu}/{abcd}, q^{1-\nu+z}/d, q^{1-\nu-z}/d
\end{array}
;\; q, \, q \r) \notag \\
&\quad+
\frac{\l(q^{-\nu}, q^{1-\nu}/{ad},  q^{1-\nu}/{bd}, q^{1-\nu}/{cd},
abcq^{\nu+z}, abcq^{\nu-z} ;q \r)_\infty}
{\l(ab, ac, bc, q^{1-2\nu}/{abcd},
q^{1+z}/d, q^{1-z}/d ;q \r)_\infty}
\notag \\
&\quad\times
\ _4\f_3\l(
\begin{array}{cc}
abcdq^{\nu-1}, abq^{\nu}, acq^{\nu}, bcq^{\nu} \\
\, \\
abcdq^{2\nu}, abcq^{\nu+z}, abcq^{\nu-z}
\end{array}
;\; q, \, q \r). \notag
\end{align}
These representations show that $u_\nu(x; a, b, c; d)$ is
actually symmetric in $a$, $b$, and $c$.
Here we have used the standard notation,
\begin{equation}
\begin{split}
&\ _{2r+2}W_{2r+1}\l(a; a_1,a_2, \ldots, a_{2r-1}; q, t\r)\\
&:=\ _{2r+2}\f_{2r+1}\l(
\begin{array}{cc}
a, q\sqrt{a}, -q\sqrt{a}, a_1, a_2, \ldots, a_{2r-1} \\
\, \\
\sqrt{a}, -\sqrt{a}, qa/a_1, qa/a_2, \ldots, qa/{a_{2r-1}}
\end{array}
;\; q, \,t\r),
\end{split}
\end{equation}
for a very-well-poised basic hypergeometric series defined before.

Let us discuss analyticity properties of the function $u_\nu(z)$
defined in (2.7)--(2.8). One can easily see, that for
integers $\nu=n=0, 1, 2, \ldots$ this function
is just a multiple of the Askey--Wilson polynomial (1.1)
of the $n$-th degree,
\begin{equation}
u_n(x; a, b, c; d)
=\frac{(-1)^n q^{-n(n-1)/2}}{(ab, ac, bc; q)_n}\,d^{-n}\;
p_n(x; a, b, c, d),
\end{equation}
where $x=\l(q^z+q^{-z}\r)/2$.
In this case the Askey--Wilson polynomials $p_n(x; a, b, c, d)$
are symmetric with respect to a permutation of all four parameters
$a$, $b$, $c$, and $d$ due to Sears' transformation \cite{As:Wi}.

On the other hand, if $\nu$ is not an integer,
the essential poles of the $_8\f_7$-solution in (2.7)
coincide with the simple poles of the infinite product
$$\l( q^{1+z}/d, q^{1-z}/d; q\r)_\infty^{-1}.$$
Therefore, the function
\begin{equation}
v_{\nu}(z)=v_{\nu}(x(z); a,b,c;d):=
{\l( q^{1+z}/d, q^{1-z}/d; q\r)_\infty}\; u_{\nu}(x(z); a, b, c; d),
\end{equation}
where $u_{\nu}(x; a, b, c; d)$ is defined in (2.7)--(2.8),
is an entire function in the complex $z$-plane.

Let us mention also
%It is not difficult to show
that function $u_{\nu}(x; a, b, c; d)$
satisfies the simple difference-differentiation formula
\begin{align}
\frac{\delta}{\delta x(z)}\, u_{\nu}(x(z); a, b, c; d)
=&\frac{2q}{(1-q)d}\,
\dfrac{\l(1-q^{-\nu}\r)\l(1-abcdq^{\nu-1}\r)}{(1-ab)(1-ac)(1-bc)} \\
&\times u_{\nu-1}\l(x(z); aq^{1/2}, bq^{1/2}, cq^{1/2}; dq^{1/2}\r),
\notag
\end{align}
where $\delta f(z)=f(z+1/2)-f(z-1/2)$ and $x(z)=\l(q^z+q^{-z}\r)/2$.

\section{Solutions of Pearson Equation}

In order to rewrite equation (2.1) for the function (2.7)--(2.8) in
self-adjoint form (2.5), we have to find a  solution of
the Pearson-type equation (2.6). In the case of the $q$-quaratic grid
$x=\dfrac{1}{2}\l(q^z+q^{-z}\r)$ this equation
can be rewritten in the form
\begin{align}
&\frac{\rho(z+1)}{\rho(z)}=
\frac{\sigma(-z)}{\sigma(z+1)} \\
&= q^{-4z-2}\, q^{2z+1}\,
\frac{(1-aq^z)(1-bq^z)(1-cq^z)(1-q^{-z}/d)}
{(1-aq^{-z-1})(1-bq^{-z-1})(1-cq^{-z-1})(1-q^{z+1}/d)}. \notag
\end{align}
It is easy to check that
\begin{align}
&\frac{\rho_0(z+1)}{\rho_0(z)}= q^{-4z-2} \qquad \text{for} \quad
\rho_0(z)=\frac{(q^{2z}, q^{-2z}; q)_\infty}{q^z-q^{-z}}, \\
&\frac{\rho_{\a}(z+1)}{\rho_{\a}(z)}= q^{-2z-1} \qquad \text{for} \quad
\rho_{\a}(z)=(\a q^z, \a q^{-z}, q^{1+z}/{\a}, q^{1-z}/{\a}; q)_\infty, \\
&\frac{\rho_a(z+1)}{\rho_a(z)}=
\frac{1-a q^{-z-1}}{1-a q^z} \qquad \text{for} \quad
\rho_a(z)=(a q^z, a q^{-z}; q)_\infty .
\end{align}
(See, for example, \cite{Ni:Su:Uv}, \cite{Ra:Su1},
and \cite{Ra:Su2} for methods of solving the Pearson equation.)
Therefore, one can choose the following solution of (3.1),
\begin{equation}
\rho(z)=
\frac{(q^z-q^{-z})^{-1}(q^{2z}, q^{-2z}, q^{1+z}/d, q^{1-z}/d; q)_\infty}
{(\a q^z, \a q^{-z}, q^{1+z}/{\a}, q^{1-z}/{\a},
aq^z, aq^{-z}, bq^z, bq^{-z}, cq^z, cq^{-z}; q)_\infty},
\end{equation}
where $\a$ is an arbitrary additional parameter.
It was shown in \cite{Su1} that this solution satisfies the correct
boundary conditions for the second order divided-difference Askey--Wilson
operator (2.5) for certain values of the parameter $\a$.

Special cases $\a=d$ or $\a=q/d$ of (3.5) result in the weight function for
the Askey--Wilson polynomials (cf. \cite{As:Wi}, \cite{At:Su1}).

One can rewrite (3.1) as
\begin{equation}
\frac{\rho(z+1)}{\rho(z)}=
\frac{(1-aq^z)(1-bq^z)(1-q^{-z}/c)(1-q^{-z}/d)}
{(1-aq^{-z-1})(1-bq^{-z-1})(1-q^{z+1}/c)(1-q^{z+1}/d)}.
\end{equation}
We shall use the corresponding solution,
\begin{equation}
\rho(z)=\frac{\l(q^{1+z}/c, q^{1-z}/c, q^{1+z}/d, q^{1-z}/d ;q\r)_\infty}
{\l(aq^z, aq^{-z}, bq^z, bq^{-z};q\r)_\infty},
\end{equation}
in Section 7.

\section{Orthogonality Property}

Now we can prove the orthogonality
relation of the $_8\f_7$-functions (2.7) with respect to
the weight function (3.5) established in \cite{Su1}.
Let us apply the following $q$-version of
the Sturm--Liouville procedure (cf. \cite{At:Su1}, \cite{Br:Ev:Is}, and
\cite{Ni:Su:Uv}).
Consider the difference equations for
the functions $u_{\nu}(z)=u_{\nu}(x(z); a, b, c; d)$ and
$u_{\mu}(z)=u_{\mu}(x(z); a, b, c; d)$ in self-adjoint form,
\begin{equation}
\frac{\Delta}{\nabla x_1(z)}\l(\sigma(z) \, \rho(z) \,
\frac{\nabla u_{\mu}(z)}{\nabla x(z)} \r)
+ \lambda_{\mu} \, \rho(z) \,u_{\mu}(z) = 0,
\end{equation}
\begin{equation}
\frac{\Delta}{\nabla x_1(z)}\l(\sigma(z) \, \rho(z) \,
\frac{\nabla u_{\nu}(z)}{\nabla x(z)} \r)
+ \lambda_{\nu} \, \rho(z) \,u_{\nu}(z) = 0,
\end{equation}
where the eigenvalues $\lambda = \lambda_{\nu}$
and $\lambda'=\lambda_{\mu}$ are defined by (2.4).
Let us multiply the first equation by $u_{\nu}(z)$,
the second one by $u_{\mu}(z)$,
and subtract the second equality from the first one.
As a result we get
\begin{equation}
\begin{split}
&\l(\lambda_{\mu} - \lambda_{\nu}\r)\;
u_{\mu}(z)\, u_{\nu}(z)\; \rho(z)\,\nabla x_1(z) \\
&\quad =\Delta \l[ \sigma(z)\, \rho(z)\; W\l(u_{\mu}(z), u_{\nu}(z)\r)\r],
\end{split}
\end{equation}
where
\begin{align}
W \l(u_{\mu}(z), \; u_{\nu}(z) \r)&=
\begin{vmatrix}
u_{\mu}(z) & u_{\nu}(z) \\
\dfrac{\nabla u_{\mu}(z)}{\nabla x(z)} &
\dfrac{\nabla u_{\nu}(z)}{\nabla x(z)}
\end{vmatrix}
\\
&=
u_{\mu}(z)\; \frac{\nabla u_{\nu}(z)}{\nabla x(z)}
-u_{\nu}(z)\; \frac{\nabla u_{\mu}(z)}{\nabla x(z)} \notag \\
&=\frac{u_{\nu}(z)\, u_{\mu}(z-1) - u_{\mu}(z)\, u_{\nu}(z-1)}
{x(z)-x(z-1)} \notag
\end{align}
is the analog of the Wronskian \cite{Ni:Su:Uv}.

We need to know the pole structure of
the analog of the Wronskian $W(u_{\mu}, u_{\nu})$
in (4.4). Let us transform the $u$'s to the entire functions $v$'s
by (2.14),
\begin{equation}
u_{\ep}(z)=\f(z)\; v_{\ep}(z),
\end{equation}
where $\ep=\mu, \nu$ and
\begin{equation}
\f(z)=\l(q^{1+z}/d, q^{1-z}/d; q\r)_\infty^{-1}.
\end{equation}
Thus,
\begin{equation}
W(u_{\mu}(z), u_{\nu}(z))= \f(z)\; \f(z-1) \;
W(v_{\mu}(z), v_{\nu}(z)),
\end{equation}
where the new ``Wronskian", $W(v_{\mu}(z), v_{\nu}(z))$, is
clearly an entire function in $z$.

Integrating (4.3) over the contour $C$ indicated in the Figure;
where the variable $z$ is such that $z=i\theta/{\log q^{-1}}$
and $-\pi \le \theta \le \pi$ (we shall assume that $0<q<1$ throughout
this work); gives
\begin{align}
\l(\lambda_{\mu} - \lambda_{\nu}\r)\;
&\int_C u_{\mu}(z)\, u_{\nu}(z)\; \rho(z)\,\nabla x_1(z)\, dz \\
= &\int_C \Delta \l[ \sigma(z)\, \rho(z)\, \f(z)\, \f(z-1)\;
W\l(v_{\mu}(z), v_{\nu}(z)\r)\r]\, dz.
\notag
\end{align}
All poles of the integrand in the right side of (4.8) coincide
with the simple poles of the function
\begin{align}
&\sigma(z)\, \rho(z)\, \f(z)\, \f(z-1) \\
&=-\frac{d\,(q^{2z}, q^{1-2z}; q)_\infty}
{\l(\a q^z, \a q^{-z}, q^{1+z}/{\a}, q^{1-z}/{\a}; q\r)_\infty} \notag\\
&\quad\times \l(aq^z, aq^{1-z},
bq^z, bq^{1-z}, cq^z, cq^{1-z},
q^{1+z}/d, q^{2-z}/d
; q\r)_\infty^{-1}.\notag
\end{align}
The integrand in the right side of (4.8) has the natural purely imaginary
period $T=2\pi i/{\log q}$ when $0<q<1$, so this integral is equal to
\begin{equation}
\int_D \l[ \sigma(z)\, \f(z)\, \f(z-1)\, \rho(z)\;
W\l(v_\mu(z), v_\nu(z)\r)\r]\, dz,
\end{equation}
where $D$ is the boundary of the rectangle on the Figure
oriented counterclockwise.

\begin{figure}

\unitlength=1.00mm
\special{em:linewidth 0.8pt}
\linethickness{0.8pt}
\begin{picture}(159.00,150.00)
\put(80,10){\line(0,1){140}}
\put(80,135){\line(1,0){55}}
\put(135,135){\line(0,-1){110}}
\put(135,25){\line(-1,0){55}}
\put(10,80){\line(1,0){140}}
\put(150,80){\line(1,0){9}}
\put(95.00,80.00){\circle*{1.50}}
\put(120.00,80.00){\circle*{1.50}}
\put(95.00,90.00){\makebox(0,0)[ct]{$\alpha_0$}}
\put(120.00,90.00){\makebox(0,0)[ct]{$1-\alpha_0$}}
\put(70.00,70.00){\makebox(0,0)[rb]{$0$}}
\put(145.00,70.00){\makebox(0,0)[lb]{$1$}}
\put(155.00,90.00){\makebox(0,0)[lt]{$\text{Re}\; z$}}
\put(90.00,150.00){\makebox(0,0)[cc]{$\text{Im}\; z$}}
\put(70.00,135.00){\makebox(0,0)[rc]{$\dfrac{i\pi}{\log q^{-1}}$}}
\put(70.00,25.00){\makebox(0,0)[rc]{$\dfrac{i\pi}{\log q}$}}
\put(145.00,135.00){\makebox(0,0)[lc]{$\dfrac{i\pi}{\log q^{-1}}+1$}}
\put(145.00,25.00){\makebox(0,0)[lc]{$\dfrac{i\pi}{\log q}+1$}}
\put(70.00,90.00){\makebox(0,0)[rc]{$C$}}
\thicklines{
\put(80.00,25.00){\vector(0,1){110.00}}
\put(80.00,135.00){\vector(0,1){0.00}}}
\end{picture}

\centerline{ \bf{Figure}}

\end{figure}

The analog of the Wronskian $W(v_\mu, v_\nu)$ is an entire function.
Thus, when $\text{max}\l(|a|, |b|, |c|, |q/d|\r)<1$,
the poles of the integrand in (4.10) inside the rectangle
in the Figure are just the simple poles of $\rho(z)$ at
$z=\alpha_0$ and $z=1-\alpha_0$, where $q^{\a_0}=\a$
and $0<\text{Re } \a_0<1/2$. By Cauchy's theorem,
\begin{align}
&\frac{1}{2\pi i}\int_D \l[ \sigma(z)\, \f(z)\, \f(z-1)\, \rho(z)\;
W\l(v_\mu(z), v_\nu(z)\r)\r]\, dz \\
&\qquad=\l.\text{Res }f(z)\r|_{z=\a_0}
+\l.\text{Res }f(z)\r|_{z=1-\a_0}, \notag
\end{align}
where
\begin{align}
f(z)=&\sigma(z)\, \rho(z)\, \f(z)\, \f(z-1) \;
W\l(v_{\mu}(z), v_{\nu}(z) \r)\\
=&-\frac{d\,(q^{2z}, q^{1-2z}; q)_\infty}
{\l(\a q^z, \a q^{-z}, q^{1+z}/{\a}, q^{1-z}/{\a}; q\r)_\infty} \notag\\
&\times
\frac{W\l(v_{\mu}(z), v_{\nu}(z)\r)}
{\l(aq^z, aq^{1-z},
bq^z, bq^{1-z}, cq^z, cq^{1-z},
q^{1+z}/d, q^{2-z}/d
; q\r)_\infty}.\notag
\end{align}
Evaluation of the residues at these simple poles gives
\begin{align}
\l.\text{Res }f(z)\r|_{z=\a_0}=&\lim_{z\to\a_0}(z-\a_0)\;f(z) \notag\\
=&-\frac{d \l. W\l(v_{\mu}(z), v_{\nu}(z)\r)\r|_{z=\a_0}}
{\log q \, \l(q, q, \a a, qa/\a, \a b, qb/\a, \a c, qc/\a,
q\a/d, q^2/{\a d} ; q \r)_\infty}
\end{align}
and
\begin{align}
\l.\text{Res }f(z)\r|_{z=1-\a_0}=&\lim_{z\to1-\a_0}(z-1+\a_0)\;f(z) \notag\\
=&-\frac{d \l. W\l(v_{\mu}(z), v_{\nu}(z)\r)\r|_{z=1-\a_0}}
{\log q \, \l(q, q, \a a, qa/\a, \a b, qb/\a, \a c, qc/\a,
q\a/d, q^2/{\a d} ; q \r)_\infty}.
\end{align}
But
\begin{equation}
\l. W\l(v_{\mu}(z), v_{\nu}(z)\r)\r|_{z=\a_0}
=\l. W\l(v_{\mu}(z), v_{\nu}(z)\r)\r|_{z=1-\a_0}
\end{equation}
due to the last line in (4.4) and
the symmetry $v_{\varepsilon}(z)=v_{\varepsilon}(-z)$,
$x(z)=x(-z)$. Thus, the residues are equal and as a result
we evaluate the integral in the right side of (4.8),
\begin{align}
&\int_C \l[ \sigma(z)\, \f(z)\, \f(z-1)\, \rho(z)\;
W\l(v_\mu(z), v_\nu(z)\r)\r]\, dz \\
&=-\frac{4\pi i\, d\; W\l(v_{\mu}(\a_0), v_{\nu}(\a_0)\r)}
{\log q \, \l(q, q, \a a, qa/\a, \a b, qb/\a, \a c, qc/\a,
q\a/d, q^2/{\a d} ; q \r)_\infty} \notag\\
&=-\frac{4\pi i\, d\; \l(\a/d, q/{\a d} ; q\r)_\infty \,
W\l(u_{\mu}(\a_0), u_{\nu}(\a_0)\r)}
{\log q \, \l(q, q, \a a, qa/\a, \a b, qb/\a, \a c, qc/\a
; q \r)_\infty} \notag
\end{align}
by (4.7).

Combining (4.8) and (4.16), we, finally, arrive at the main equation,
\begin{align}
&\l(\lambda_\mu-\lambda_\nu\r)\;
\int_C u_\mu(z)\, u_\nu(z)\; \rho(z)\,\delta x(z)\, dz \\
&=-\frac{4 \pi i \, d\; \l(\a/d, q/{\a d};q\r)_\infty}
{\log q \; (q, q, \a a, q a/\a, \a b, qb/\a, \a c, qc/\a; q)_\infty} \notag\\
&\qquad\times
W\l(u_\mu(\alpha_0), u_\nu(\alpha_0)\r), \notag
\end{align}
where $\text{max}\l(|a|, |b|, |c|, |q/d|\r)<1$
and $0< \text{Re } \a_0 < 1/2$, $q^{\a_0}=\a$.

It is worth mentioning the important special case first.
If both of the ``degree'' parameters $\mu$ and $\nu$
are nonnegative integers:
$\mu=m=0, 1, 2, \ldots$ and $\nu=n=0, 1, 2, \ldots$;
equations (4.17) and (2.13) result in the following real integral
\begin{align}
&\int_0^{\pi}
\frac{p_m(\cos \t ; a, b, c; d) \; p_n(\cos \t ; a, b, c; d)}
{\l(a e^{i\t},  a e^{-i\t}, b e^{i\t}, b e^{-i\t}, c e^{i\t},  c e^{-i\t}
; q\r)_\infty}\\
&\qquad\times
\frac{\l(e^{2i\t}, \, e^{-2i\t}, q e^{i\t}/d,  q e^{-i\t}/d
; q\r)_\infty}
{\l(\a e^{i\t}, \a e^{-i\t}, q e^{i\t}/\a, q e^{-i\t}/\a; q\r)_\infty}
\, d \t \notag \\
&\quad=
\frac{\l(\a/d, q/\a d; q\r)_\infty}
{\l(q, q, \a a, qa/\a, \a b, qb/\a , \a c, qc/\a; q\r)_\infty}
\notag \\
&\qquad\times
\frac{-4\pi q^{1/2} d}{1-q} \;\;
\frac
{W\l(p_m\l(\eta; a, b, c; d\r),\;
p_n\l(\eta; a, b, c; d\r)\r) }
{\lambda_m -\lambda_n} \notag
\end{align}
involving the Askey--Wilson polynomials.
Here we use the notation
\begin{equation}
\eta=x(\a_0)=\dfrac{1}{2}\l(\a+\a^{-1}\r).
\end{equation}
One can easily see that
when $\a=d$ or $\a=q/d$ our equation (4.18) implies the orthogonality
relation for the Askey--Wilson polynomials (1.7).

The new orthogonality property for the $_8\phi_7$-function appears
when both of the parameters $\mu$ and $\nu$ are not nonnegative integers.
It is convenient in this case to rewrite (4.17) in terms of
the entire functions $v$'s as follows
\begin{align}
&\int_0^{\pi}
\frac{v_{\mu}(\cos \t ; a, b, c; d) \; v_{\nu}(\cos \t ; a, b, c; d)}
{\l(a e^{i\t},  a e^{-i\t}, b e^{i\t}, b e^{-i\t}, c e^{i\t},  c e^{-i\t},
q e^{i\t}/d,  q e^{-i\t}/d; q\r)_\infty}\\
&\qquad\times
\frac{\l(e^{2i\t}, \, e^{-2i\t}; q\r)_\infty}
{\l(\a e^{i\t}, \a e^{-i\t}, q e^{i\t}/\a, q e^{-i\t}/\a; q\r)_\infty}
\, d \t \notag \\
&\quad=\l(q, q, \a a, qa/\a, \a b, qb/\a , \a c, qc/\a,
q\a/d, q^2/{\a d}; q\r)_\infty^{-1} \notag \\
&\qquad\times
\frac{-4\pi q^{1/2} d}{1-q} \;\;
\frac
{W\l(v_\mu\l(\eta; a, b, c; d\r),\;
v_\nu\l(\eta; a, b, c; d\r)\r) }
{\lambda_\mu -\lambda_\nu}. \notag
\end{align}
The limiting case $\mu \to \nu$ of (4.20) is also of interest.
From (4.4),
\begin{align}
&\lim_{\mu\to\nu}\;
\frac
{W\l(v_\mu(z),\; v_\nu(z)\r)}
{\lambda_\mu -\lambda_\nu} \notag\\
&=\lim_{\mu\to\nu}\;
\frac
{v_\mu(z)\, \dfrac{\nabla v_\nu(z)}{\nabla x(z)}
-v_\nu(z)\, \dfrac{\nabla v_\mu(z)}{\nabla x(z)}}
{\lambda_\mu -\lambda_\nu} \notag\\
&=\lim_{\mu\to\nu}\;
\frac
{\dfrac{\partial v_\mu(z)}{\partial \mu} \,
\dfrac{\nabla v_\nu(z)}{\nabla x(z)}
-v_\nu(z)\, \dfrac{\partial}{\partial \mu}
\l(\dfrac{\nabla v_\mu(z)}{\nabla x(z)}\r)}
{\partial\lambda_\mu/{\partial\mu}} \notag\\
&=
\frac
{\dfrac{\partial v_\nu(z)}{\partial \nu} \,
\dfrac{\nabla v_\nu(z)}{\nabla x(z)}
-v_\nu(z)\, \dfrac{\partial}{\partial \nu}
\l(\dfrac{\nabla v_\nu(z)}{\nabla x(z)}\r)}
{\partial\lambda_\nu/{\partial\nu}}. \notag
\end{align}
Therefore,
\begin{align}
&\int_0^{\pi}
\frac{\l(v_{\nu}(\cos \t ; a, b, c; d)\r)^2}
{\l(a e^{i\t},  a e^{-i\t}, b e^{i\t}, b e^{-i\t}, c e^{i\t},  c e^{-i\t},
q e^{i\t}/d,  q e^{-i\t}/d; q\r)_\infty}\\
&\qquad\times
\frac{\l(e^{2i\t}, \, e^{-2i\t}; q\r)_\infty}
{\l(\a e^{i\t}, \a e^{-i\t}, q e^{i\t}/\a, q e^{-i\t}/\a; q\r)_\infty}
\, d \t \notag \\
&\quad=\l(q, q, \a a, qa/\a, \a b, qb/\a , \a c, qc/\a,
q\a/d, q^2/{\a d}; q\r)_\infty^{-1} \notag \\
&\qquad\times
\frac{-4\pi q^{1/2} d}{1-q} \;
\Bigg[\frac{\partial}{\partial \lambda_\nu}\,
v_\nu(\eta; a, b, c; d)\;
\l(\frac{\nabla}{\nabla x}\,
v_\nu(x; a, b, c; d)\r)\Bigg|_{x=\eta} \notag \\
&\qquad\qquad- v_\nu(\eta; a, b, c; d)\;
\frac{\partial}{\partial \lambda_\nu}
\l(\frac{\nabla}{\nabla x}\,
v_\nu(x; a, b, c; d)\r) \Bigg|_{x=\eta}\Bigg]. \notag
\end{align}

Finally, choosing the parameters $\mu$ and $\nu$ as $\ep$-solutions
of the equation
\begin{equation}
v_{\ep}\l(\dfrac{1}{2}\l(\a+\a^{-1}\r); a, b, c; d\r)=0,
\end{equation}
we arrive from (4.20) and (4.21) at the {\it orthogonality relation\/} of
the $_8\f_7$-functions under consideration,
\begin{align}
&\int_0^{\pi}
\frac{v_{\mu}(\cos \t ; a, b, c; d) \;
v_{\nu}(\cos \t ; a, b, c; d)}
{\l(a e^{i\t},  a e^{-i\t}, b e^{i\t}, b e^{-i\t}, c e^{i\t},  c e^{-i\t},
q e^{i\t}/d,  q e^{-i\t}/d; q\r)_\infty}\\
&\qquad\times
\frac{\l(e^{2i\t}, \, e^{-2i\t}; q\r)_\infty}
{\l(\a e^{i\t}, \a e^{-i\t}, q e^{i\t}/\a, q e^{-i\t}/\a; q\r)_\infty}
\, d \t =0 \notag
\end{align}
if $\mu \ne \nu$, and
\begin{align}
&\int_0^{\pi}
\frac{\l(v_{\nu}(\cos \t ; a, b, c; d)\r)^2}
{\l(a e^{i\t},  a e^{-i\t}, b e^{i\t}, b e^{-i\t}, c e^{i\t},  c e^{-i\t},
q e^{i\t}/d,  q e^{-i\t}/d; q\r)_\infty}\\
&\qquad\times
\frac{\l(e^{2i\t}, \, e^{-2i\t}; q\r)_\infty}
{\l(\a e^{i\t}, \a e^{-i\t}, q e^{i\t}/\a, q e^{-i\t}/\a; q\r)_\infty}
\, d \t \notag \\
&\quad=\l(q, q, \a a, qa/\a, \a b, qb/\a , \a c, qc/\a,
q\a/d, q^2/{\a d}; q\r)_\infty^{-1} \notag \\
&\qquad\times
\frac{-4\pi q^{1/2} d}{1-q} \;\;
\frac{\partial}{\partial \lambda_\nu}\,
v_\nu(\eta; a, b, c; d)\;
\l(\frac{\nabla}{\nabla x}\,
v_\nu(x; a, b, c; d)\r)\Bigg|_{x=\eta} \notag
\end{align}
if $\mu=\nu$, respectively.
We remind the reader that
$\text{max}\l(|a|, |b|, |c|, |q/d|\r)<1$.
This will be assumed throughout this work.

One can easily see that $v_{\varepsilon}(x; a, b, c; d)$ with
$\varepsilon=\mu, \nu$ in (4.23)--(4.24) are
real-valued functions of $x$ for real $\varepsilon$
which are orthogonal with respect to a positive weight function
when all parameters $a$, $b$, $c$, $d$, and $\a$ are real,
or when any two of $a$, $b$, $c$ are complex conjugate and all other
parameters are real.
We shall usually assume that as well.

\section{Properties of Zeros and Asymptotics}

In the previous section we have proved that
the $_8\f_7$-functions (2.7)--(2.8) are orthogonal if
the ``boundary'' condition (4.22) is satisfied.
Let us discuss now some properties of $\nu$-zeros of
the corresponding ``boundary'' function,
\begin{align}
&v_{\nu}\l(\frac{1}{2}\l(\a+\a^{-1}\r);a, b, c; d\r) \notag \\
&=\frac
{\l(qa/d, bcq^\nu, \a q^{1-\nu}/d, q^{1-\nu}/{\a d}; q\r)_\infty}
{\l(q^{1-\nu}a/d, bc; q\r)_\infty} \\
&\quad\times \ _8\f_7\l(
\begin{array}{cc}
\dfrac{aq^{-\nu}}{d},
q \sqrt{\dfrac{aq^{-\nu}}{d}},-q \sqrt{\dfrac{aq^{-\nu}}{d}},
q^{-\nu}, \dfrac{q^{1-\nu}}{bd}, \dfrac{q^{1-\nu}}{cd}, \a a, a/\a \\
\, \\
\sqrt{\dfrac{aq^{-\nu}}{d}}, -\sqrt{\dfrac{aq^{-\nu}}{d}},
\dfrac{aq}{d}, ab, ac, \dfrac{q^{1-\nu}\a}{d}, \dfrac{q^{1-\nu}}{\a d}
\end{array}
;\; q, \, bcq^{\nu} \r) \notag\\
&=\frac
{\l(q\a/d, q^{1-\nu}/{\a d}, abcq^{\nu}/\a; q\r)_\infty}
{\l(abc/{\a}; q\r)_\infty} \\
&\quad\times \ _8\f_7\l(
\begin{array}{cc}
\dfrac{abc}{q \a},
\sqrt{\dfrac{abc}{q^{1/2} \a}}, -\sqrt{\dfrac{abc}{q^{1/2} \a}},
q^{-\nu}, abcd q^{\nu-1}, \dfrac{a}{\a}, \dfrac{b}{\a}, \dfrac{c}{\a} \\
\, \\
\sqrt{\dfrac{abc}{q \a}}, -\sqrt{\dfrac{abc}{q \a}},
\dfrac{q^{\nu}abc}{\a}, \dfrac{q^{1-\nu}}{\a d}, ab, ac, bc
\end{array}
;\; q, \, q\,\dfrac{\a}{d} \r). \notag
\end{align}
We have used (III.23) of \cite{Ga:Ra} to transform (5.1) to (5.2),
compare also (2.10).
Again, $v_\nu(\eta; a, b, c; d)$ is real-valued function
of $\nu$ when all parameters $a$, $b$, $c$, $d$, and $\a$ are real,
or when any two of $a$, $b$, $c$ are complex conjugate and all other
parameters are real.

Main properties of zeros of the function (5.1)--(5.2) can be investigated
by using the same methods as in \cite{Bu:Su}, \cite{Is}, \cite{Is:Ma:Su2}
and \cite{Is:Ma:Su3}.
The first property is that under certain conditions the real-valued function
$v_{\nu}(\eta; a, b, c; d)$
has an infinity of positive $\nu$-zeros.
In order to establish this fact, one can consider
the large $\nu$-asymptotics of the $_8\f_7$-function in (5.2),
\begin{align}
&\ _8\f_7\l(
\begin{array}{cc}
\dfrac{abc}{q \a},
\sqrt{\dfrac{abc}{q^{1/2} \a}}, -\sqrt{\dfrac{abc}{q^{1/2} \a}},
q^{-\nu}, abcd q^{\nu-1}, \dfrac{a}{\a}, \dfrac{b}{\a}, \dfrac{c}{\a} \\
\, \\
\sqrt{\dfrac{abc}{q \a}}, -\sqrt{\dfrac{abc}{q \a}},
\dfrac{q^{\nu}abc}{\a}, \dfrac{q^{1-\nu}}{\a d}, ab, ac, bc
\end{array}
;\; q, \, q\,\dfrac{\a}{d} \r) \\
&\qquad@>>{ \nu \to \infty}>
\ _6\f_5\l(
\begin{array}{cc}
\dfrac{abc}{q \a},
\sqrt{\dfrac{abc}{q^{1/2} \a}}, -\sqrt{\dfrac{abc}{q^{1/2} \a}},
\dfrac{a}{\a}, \dfrac{b}{\a}, \dfrac{c}{\a} \\
\, \\
\sqrt{\dfrac{abc}{q \a}}, -\sqrt{\dfrac{abc}{q \a}},
ab, ac, bc
\end{array}
;\; q, \,{\a}^2 \r) \notag\\
&\qquad\qquad\qquad =
\frac{\l(\a a, \a b, \a c, abc/\a; q\r)_\infty}
{\l(ab, ac, bc, \a^2; q \r)_\infty}
\notag
\end{align}
by (II.20) of \cite{Ga:Ra}. Therefore,
\begin{align}
v_{\nu}\l(\frac{1}{2}\l(\a+\a^{-1}\r); a, b, c; d\r)
=&\frac{\l(\a a, \a b, \a c,  q\a/d; q\r)_\infty}
{\l(ab, ac, bc, \a^2; q\r)_\infty} \\
&\times \l(q^{1-\nu}/{\a d} ; q \r)_\infty \,
[1+\text{o}(1)], \notag
\end{align}
as $\nu\to\infty$.
But for the positive values of $q$ and $\a d$ the function
\begin{equation*}
\l(q^{1-\nu}/{\a d} ; q \r)_\infty
\end{equation*}
oscillates and has an infinity of real zeros as $\nu$ approaches infinity
(see \cite{Is} and \cite{Is:Ma:Su2}).
Indeed, consider the points $\nu=\omega_n={\omega}_0+n$,
such that $q^{{\omega}_0}=\b$,
where $n=0, 1, 2, \ldots $ and $q<\b<1$, as test points.
Then, by using (I.8) of \cite{Ga:Ra},
\begin{align}
v_{{\omega}_0+n}\l(\eta; a, b, c; d\r)
=&\frac{\l(\a a, \a b, \a c,  q\a/d; q\r)_\infty}
{\l(ab, ac, bc, \a^2; q\r)_\infty} \\
&\times (-1)^n\, q^{-n(n+1)/2}\;
\frac{(\a \b dq ; q)_n}{(\a \b d)^n}\;
[1+\text{o}(1)], \notag
\end{align}
as $n\to\infty$, and one can see that the right side of (5.5) changes sign
infinitely many times at the test points $\nu={\omega}_n$ as $\nu$ approaches
infinity.

Thus we have established the following theorem.

\begin{theorem}
The real-valued function $v_{\nu}\l(\eta; a, b, c; d\r)$ defined by (5.2)
has an infinity of positive $\nu$-zeros when
$\a q<d$, $\a^2<1$, $\a d>0$, and $0<q<1$.
Also, all parameters $a$, $b$, $c$ are real, or when any of
two of $a$, $b$, $c$ are complex conjugate and the remainder parameter
is real.
\end{theorem}

Now we can prove the next result.

\begin{theorem}
Function $v_{\nu}\l(\eta; a, b, c; d\r)$ defined in (5.2)
has only real $\nu$-zeros under the following conditions:

\begin{enumerate}

\item parameters $\a$ and $d$ are real; parameters $a$, $b$, $c$ are all
      real, or, if any two of them are complex conjugate, the third
      parameter is real;

\item the following inequalities holds:

\begin{enumerate}

\item $\text{max}(|a|, |b|, |c|, |q/d|)<1, \quad q^{1/2}<\a<1;$

\item $ \a d>q, \quad \a abc<q.$

\end{enumerate}

\end{enumerate}

\end{theorem}

\begin{pf}
Suppose that $\nu_0$ is a zero of function (5.2) which is not real.
Let $\nu_1$ be the complex number conjugate to $\nu_0$,
so that $\nu_1$ is also a zero of (5.2) because this function is real
under the hypotheses of the theorem.

Consider equation (4.20) with $\mu=\nu_0$ and $\nu=\nu_1$.
The integral on the left does not equal zero due to the positivity
of the integrand, but the analog of the Wronskian on the right is zero
in view of (4.22). Therefore,
\begin{equation}
\lambda_{\nu_0}=\lambda_{\nu_1},
\end{equation}
the eigenvalues defined by (2.4) are real.
The last equation can be rewritten as
\begin{equation}
\l(1-q^{\nu_1-\nu_0}\r)\; \l(1-abcd q^{\nu_0+\nu_1-1}\r)=0.
\end{equation}

The first solution is $\nu_0=\nu_1$, so the roots are real in this case.

The second solution of (5.7) is
\begin{equation}
q^{\nu_0}=\sqrt{\frac{q}{abcd}}\; e^{i\chi}, \qquad
q^{\nu_1}=\sqrt{\frac{q}{abcd}}\; e^{-i\chi},
\end{equation}
where $abc>0$ and $\chi$ is an arbitrary real number.
But in this case our function can be represented as a multiple of
a positive function. Indeed, when $\nu=\nu_0$ or $\nu=\nu_1$,
\begin{align}
&v_{\nu}(\eta; a, b, c; d) \\
&=\frac
{\l(q/\a d, q^{1-\nu}\a/d, \a abcq^{\nu}; q\r)_\infty}
{\l(\a abc; q\r)_\infty} \notag \\
&\quad\times \ _8\f_7\l(
\begin{array}{cc}
\dfrac{\a abc}{q},
\sqrt{\dfrac{\a abc}{q^{1/2}}}, -\sqrt{\dfrac{\a abc}{q^{1/2}}},
q^{-\nu}, abcd q^{\nu-1}, \a a, \a b, \a c \\
\, \\
\sqrt{\dfrac{\a abc}{q}}, -\sqrt{\dfrac{\a abc}{q}},
q^{\nu} \a abc, q^{1-\nu} \a/d, ab, ac, bc
\end{array}
;\; q, \, \dfrac{q}{\a d} \r) \notag\\
&=\frac
{\l(q/\a d, q^{1-\nu}\a/d, \a abcq^{\nu}; q\r)_\infty}
{\l(\a abc; q\r)_\infty} \notag \\
&\quad\times\Bigg[
1+\sum_{n=1}^\infty
\frac{\l(1-\a abcq^{2n-1}\r) \l(\a abcq^{-1}, \a a, \a b, \a c; q\r)_n}
{\l(1-\a abcq^{-1}\r) \l(q, ab, ac, bc; q\r)_n} \notag \\
&\qquad\qquad\qquad\qquad\times\frac
{\l(q^{-\nu}, abcdq^{\nu-1}; q\r)_n}
{\l(q^{\nu} \a abc, q^{1-\nu} \a/d; q\r)_n}\;
\l(\dfrac{q}{\a d}\r)^n \Bigg] \notag
\end{align}
by (5.2), (III.23) and (I.31) of \cite{Ga:Ra}.
One can easily see that under the hypotheses of the theorem all products
in the last but one line of (5.9) are positive.
From (1.3)--(1.4) and (5.8),
\begin{align}
\l(q^{-\nu}, abcdq^{\nu-1}; q\r)_n
&=\prod_{k=0}^{n-1}
\l(1-2 \sqrt{\frac{abcd}{q}} \cos \chi \; q^k
+\frac{abcd}{q}\; q^{2k} \r) \notag \\
&=\prod_{k=0}^{n-1}
\l(\l(1-\sqrt{\frac{abcd}{q}} \cos \chi \; q^k \r)^2
+\frac{abcd}{q} \sin^2 \chi \; q^{2k} \r)
\end{align}
and
\begin{align}
\l(q^{\nu}\a abc, q^{1-\nu}\a/d; q\r)_n
&=\prod_{k=0}^{n-1}
\l(1-2\a \sqrt{\frac{abcq}{d}} \cos \chi \; q^k
+\a^2 \frac{abcq}{d} \; q^{2k} \r) \notag \\
&=\prod_{k=0}^{n-1}
\l(\l(1-\a \sqrt{\frac{abcq}{d}} \cos \chi \; q^k \r)^2
+\a^2 \frac{abcq}{d} \sin^2 \chi \; q^{2k} \r),
\end{align}
where $n=1, 2, 3, \ldots, \infty$.
Thus, these products are positive too, which proves the positivity of
the ``boundary'' function in (5.9). So we have obtained a contradiction
and the complex zeros (5.8) cannot exist.
This completes the proof of the theorem.
\end{pf}

\begin{theorem}
The real $\nu$-zeros of $v_{\nu}(\eta; a, b, c; d)$ are simple
under the hypotheses 1 and 2(a) of Theorem 5.2.
\end{theorem}

\begin{pf}
This follows directly from the relation (4.21).
Indeed, the integral on the left side is positive under the conditions
of the theorem, which means that
$\dfrac{\partial}{\partial \nu}\, v_\nu(\eta; a, b, c; d)\ne 0$
when $v_\nu(\eta; a, b, c; d)=0$.
\end{pf}

Let us consider two functions,
\begin{align}
&f(\nu)=v_\nu(\eta; a, b, c; d), \\
&g(\nu)=\l. \frac{\nabla}{\nabla x}\, v_\nu(x; a, b, c; d)\r|_{x=\eta}.
\end{align}
Our next property is that positive zeros of $f(\nu)$ and $g(\nu)$
are interlaced.

\begin{theorem}
If $\nu_1, \nu_2, \nu_3, \ldots$ are the positive zeros of $f(\nu)$
arranged in ascending order of magnitude, and
$\mu_1, \mu_2, \mu_3, \ldots$ are those of $g(\nu)$,
then
\begin{equation}
0<\mu_1<\nu_1<\mu_2<\nu_2<\mu_3<\nu_3< \ldots,
\end{equation}
if the conditions of Theorem 5.3 are satisfied.
\end{theorem}

\begin{pf}
Suppose that $\nu_k$ and $\nu_{k+1}$ are two successive zeros of $f(\nu)$.
Then the derivative $\dfrac{\partial}{\partial \nu}f(\nu)$ has different
signs at $\nu=\nu_k$ and $\nu=\nu_{k+1}$.
This means, in view of the positivity of the integral on the left side
of (4.21), that $g(\nu)$ changes its sign between
$\nu_k$ and $\nu_{k+1}$ and, therefore, has at least one zero
on each interval $(\nu_k, \nu_{k+1})$.

To complete the proof of the theorem, we have to show that $g(\nu)$
changes its sign on each interval $(\nu_k, \nu_{k+1})$ only once.
Suppose that $g(\mu_k)=g(\mu_{k+1})=0$ and
$\nu_k<\mu_k<\mu_{k+1}<\nu_{k+1}$.
Then, by (4.21), function $f(\nu)$ has different signs at
$\nu=\mu_k$ and $\nu=\mu_{k+1}$ and, therefore, this function has
at least one more zero on $(\nu_k, \nu_{k+1})$. So, we have obtained
a contradiction, and, therefore,  function $g(\nu)$ has exactly one zero
between any two successive zeros of $f(\nu)$.

In order to finish the proof of our theorem one has to show that
$\mu_1<\nu_1$. We would like to leave the proof of this fact
as a congecture for the reader.
\end{pf}

Let us discuss an asymptotic behavior of $u_\nu(x; a, b, c; d)$
for large values of $\nu$ applying the methods used
in \cite{Is:Wi}, \cite{Ra0}, and \cite{Ga:Ra}
at the level of the Askey--Wilson polynomials.
In a similar manner, the $_8\f_7$-function in (2.10) can be transformed
by (III.37) of \cite{Ga:Ra},
\begin{align}
&\ _8W_7\l(abcq^{z-1}; q^{-\nu}, abcdq^{\nu-1}, aq^z, bq^z, cq^z;
q,\; \frac{q^{1-z}}{d}\r) \\
&=\frac{\l(aq^{-z}, bq^{-z}, cq^{-z}, bq^{1+\nu+z}, cq^{1+\nu+z},
dq^{\nu-z}, abcq^z, bcdq^{\nu+z}; q\r)_\infty}
{\l(q^{-2z}, ab, ac, bc, q^{1+\nu}, bdq^\nu, cdq^\nu, bcq^{1+\nu+2z}
; q\r)_\infty} \notag \\
&\quad\times\ _8W_7\l(bcq^{\nu+2z}; bcq^\nu, bq^z, cq^z,
q^{1+z}/a, q^{1+z}/d; q,\; adq^\nu  \r) \notag \\
&+\frac{\l(aq^z, bq^z, cq^z, bq^{1+\nu-z}, cq^{1+\nu-z},
dq^{\nu+z}, abcq^z, bcdq^{\nu-z}; q\r)_\infty}
{\l(q^{2z}, ab, ac, bc, q^{1+\nu}, bdq^\nu, cdq^\nu, bcq^{1+\nu-2z}
; q\r)_\infty} \notag \\
&\quad\times\ \frac{\l(q^{1-\nu-z}/d, q^{1+z}/d, abcq^{\nu-z}
; q\r)_\infty}
{\l(q^{1-z}/d, q^{1-\nu+z}/d, abcq^{\nu+z}; q\r)_\infty} \notag \\
&\qquad\times\ _8W_7\l(bcq^{\nu-2z}; bcq^\nu, bq^{-z}, cq^{-z},
q^{1-z}/a, q^{1-z}/d; q,\; adq^\nu  \r). \notag
\end{align}
Therefore,
\begin{align}
&u_\nu(x; a, b, c; d) \\
&=\frac{\l(bq^{1+\nu+z}, cq^{1+\nu+z}, dq^{\nu-z},
abcq^{\nu+z}, bcdq^{\nu+z}; q\r)_\infty}
{\l(q^{1+\nu}, bdq^\nu, cdq^\nu, bcq^{1+\nu+2z}
; q\r)_\infty} \notag \\
&\quad\ \times
\frac{\l(aq^{-z}, bq^{-z}, cq^{-z}, q^{1-\nu+z}/d; q\r)_\infty}
{\l(ab, ac, bc, q^{-2z}, q^{1+z}/d; q\r)_\infty} \notag \\
&\qquad\times\ _8W_7\l(bcq^{\nu+2z}; bcq^\nu, bq^z, cq^z,
q^{1+z}/a, q^{1+z}/d; q,\; adq^\nu  \r) \notag \\
&+\frac{\l(bq^{1+\nu-z}, cq^{1+\nu-z}, dq^{\nu+z},
abcq^{\nu-z}, bcdq^{\nu-z}; q\r)_\infty}
{\l(q^{1+\nu}, bdq^\nu, cdq^\nu, bcq^{1+\nu-2z}
; q\r)_\infty} \notag \\
&\quad\ \times
\frac{\l(aq^z, bq^z, cq^z, q^{1-\nu-z}/d; q\r)_\infty}
{\l(ab, ac, bc, q^{2z}, q^{1-z}/d; q\r)_\infty} \notag \\
&\qquad\times\ _8W_7\l(bcq^{\nu-2z}; bcq^\nu, bq^{-z}, cq^{-z},
q^{1-z}/a, q^{1-z}/d; q,\; adq^\nu \r). \notag
\end{align}
We shall use the last expression to determine
the large $\nu$ asymptotic of our orthogonal $_8\f_7$-function later.

For the ``boundary" function (5.2) equation (5.16) takes the form
\begin{align}
&v_\nu\l(\frac{1}{2}\l(\a+\a^{-1}\r); a, b, c; d\r) \\
&=\frac{\l(\a bq^{1+\nu}, \a cq^{1+\nu}, dq^{\nu}/\a,
\a abcq^{\nu}, \a bcdq^{\nu}; q\r)_\infty}
{\l(q^{1+\nu}, bdq^\nu, cdq^\nu, \a^2 bcq^{1+\nu}
; q\r)_\infty} \notag \\
&\quad\ \times
\frac{\l(a/\a, b/\a, c/\a, q/\a d, \a q^{1-\nu}/d; q\r)_\infty}
{\l(ab, ac, bc, \a^{-2}; q\r)_\infty} \notag \\
&\qquad\times\ _8W_7\l(\a^2 bcq^{\nu}; bcq^\nu, \a b, \a c,
q\a/a, q\a/d; q,\; adq^\nu  \r) \notag \\
&+\frac{\l(bq^{1+\nu}/\a, cq^{1+\nu}/\a, \a dq^{\nu},
abcq^{\nu}/\a, bcdq^{\nu}/\a; q\r)_\infty}
{\l(q^{1+\nu}, bdq^\nu, cdq^\nu, bcq^{1+\nu}/\a^2
; q\r)_\infty} \notag \\
&\quad\ \times
\frac{\l(\a a, \a b, \a c, q\a/d, q^{1-\nu}/\a d; q\r)_\infty}
{\l(ab, ac, bc, \a^2; q\r)_\infty} \notag \\
&\qquad\times\ _8W_7\l(bcq^{\nu}/\a^2; bcq^\nu, b/\a, c/\a,
q/\a a, q/\a d; q,\; adq^\nu \r). \notag
\end{align}
The last term dominates here as $\nu$ approaches infinity (cf. (5.4)).

The proof of Theorem 5.1 has strongly indicated that asymptotically
the large positive $\nu$-zeros of $v_\nu(\eta; a, b, c; d)$ are
\begin{equation}
\nu_n=n+{\epsilon}_n, \qquad 0\le{\epsilon}_n<1
\end{equation}
as $n\to\infty$. The same consideration as in \cite{Is}, \cite{Is:Ma:Su3},
and \cite{Bu:Su} shows that this function changes sign  only once between
any two successive test points $\nu=\omega_n$ and $\nu=\omega_{n+1}$
defined on the page 15 for sufficiently large values of $n$.
We include details of this proof in Section 9 to make this work
as self-contained as possible.

Our next theorem provides a more accurate estimate for the distribution of
the large positive zeros of the ``boundary" function (5.2).

\begin{theorem}
If $\nu_1, \nu_2, \nu_3, \ldots$ are the positive zeros of
$v_\nu(\eta; a, b, c; d)$ arranged in ascending order of magnitude,
then
\begin{equation}
\nu_n=n-\frac{\log (\a d)}{\log q}+\text{o}(1),
\end{equation}
as $n\to\infty$.
\end{theorem}

\begin{pf}
All zeros of our function (5.2) coincide with the zeros of a new function,
\begin{equation}
w_\nu(\eta; a, b, c; d):=
\l(-q^{1-\nu}/\a d; q\r)_\infty^{-1}\;
v_\nu(\eta; a, b, c; d),
\end{equation}
because $\l(-q^{1-\nu}/\a d; q\r)_\infty$ is positive for real $\nu$.
Equations (5.17) and (5.20) give us the large $\nu$-asymptotic of
$w_\nu(\eta; a, b, c; d)$.
When $\nu=\g_n$ such that $q^{\g_n}=q^n/\a d$ and $n=1, 2, 3, \ldots$,
the second term in (5.17) vanishes and we get
\begin{align}
w_{\g_n}(\eta; a, b, c; d)
&=\frac
{\l(q\a^2, a/\a, bc, qb/d, qc/d,abc/d ; q\r)_\infty}
{\l(-q, ab, ac, bc, \a bcq/d; q\r)_\infty} \\
&\times\ \frac
{\l(q/\a d, b/\a, c/\a, \a bcq/d ; q\r)_n}
{\l(-1, bc, b/d, qc/d, abc/d; q\r)_n}\;
\l(-\a^2\r)^n \notag \\
&\times\ _8W_7\l(\frac{\a bc}{d}\,q^n; \frac{b}{\a}\,q^n,
\a b, \a c, q\,\frac{\a}{a}, q\,\frac{\a}{d}; q,\; \frac{a}{\a}\,q^n\r)
\notag
\end{align}
with the help of (I.9) of \cite{Ga:Ra}.
Thus,
\begin{equation}
\lim_{n\to\infty} w_{\g_n}(\eta; a, b, c; d)=0,
\end{equation}
which proves our theorem.
\end{pf}

Let us discuss also the large positive $\nu$-asymptotic of
$u_\nu(x; a, b, c; d)$ when $x=\cos \t$ belongs to the interval of
orthogonality $-1<x<1$. It is clear from (5.16) that the leading terms
in the asymptotic expansion of this function are given by
\begin{align}
u_\nu(\cos \t; a, b, c; d)
&\sim\frac{\l(ae^{i\t}, be^{i\t}, ce^{i\t}, q^{1-\nu}e^{-i\t}/d; q\r)_\infty}
{\l(ab, ac, bc, e^{2i\t}, qe^{-i\t}/d; q\r)_\infty}\\
&+\frac{\l(ae^{-i\t}, be^{-i\t}, ce^{-i\t}, q^{1-\nu}e^{i\t}/d; q\r)_\infty}
{\l(ab, ac, bc, e^{-2i\t}, qe^{i\t}/d; q\r)_\infty}. \notag
\end{align}
In particular, when $\nu=\nu_n$ are large zeros of
$u_\nu(\eta; a, b, c; d)$,
we can estimate
\begin{equation}
u_{\nu_n}(\cos \t; a, b, c; d)
\sim
u_{\g_n}(\cos \t; a, b, c; d) ,
\end{equation}
where $q^{\g_n}=q^n/\a d$, due to (5.19) as $n\to\infty$.
Relations (5.23)--(5.24) lead to the following theorem.

\begin{theorem}
For $-1<x=\cos \t<1$ and $|q|<1$ the leading term in
the asymptotic expansion of $u_{\g_n}(\cos \t; a, b, c; d)$
as $n\to\infty$ is given by
\begin{align}
u_{\g_n}(\cos \t; a, b, c; d)
&\sim\frac{q^{-(n-1)(n-2)/2}}{\l(ab, ac, bc;q\r)_\infty}\\
&\times \l|A\l(e^{i\t}\r)\r|\;
\cos\l((n-1)\t-\chi\r), \notag
\end{align}
where
\begin{equation}
A\l(e^{i\t}\r)=\frac
{\l(ae^{i\t}, be^{i\t}, ce^{i\t}, \a e^{-i\t}, qe^{i\t}/\a; q\r)_\infty}
{\l(e^{2i\t}, qe^{-i\t}/d; q\r)_\infty},
\end{equation}
\begin{align}
\l|A\l(e^{i\t}\r)\r|^{-2}
=&\frac{\l(e^{2i\t},\; e^{-2i\t}; q\r)_\infty}
{\l(a e^{i\t}, a e^{-i\t}, b e^{i\t}, b e^{-i\t},
ce^{i\t}, ce^{-i\t}; q\r)_\infty} \\
&\times
\frac{\l(qe^{i\t}/d, \; qe^{-i\t}/d; q\r)_\infty}
{\l(\a e^{i\t}, \a e^{-i\t}, q e^{i\t}/\a, q e^{-i\t}/\a; q\r)_\infty},
\notag
\end{align}
and
\begin{equation}
\chi=\arg A\l(e^{i\t}\r).
\end{equation}
\end{theorem}

It is worth mentioning that (5.27) coinsides with the weight function
in the orthogonality relation for $u_\nu(\cos \t; a, b, c; d)$.

In a similar fashion, one can consider some other properties of zeros of
the $_8\varphi_7$-function (5.1)--(5.2)
close to those established in \cite{Bu:Su}, \cite{Is:Ma:Su2} and
\cite{Is:Ma:Su3} at the level of the basic trigonometric functions and
the $q$-Bessel function, respectively.

\section{Evaluation of Some Constants}

In this section we shall find an explicit expression for
the squared norm,
\begin{align}
d^2_\nu=&\int_0^{\pi}
\frac{u_{\nu}^2(\cos \t ; a, b, c; d)}
{\l(a e^{i\t}, a e^{-i\t}, b e^{i\t}, b e^{-i\t}, c e^{i\t},  c e^{-i\t}
; q\r)_\infty}\\
&\quad\times
\frac{\l(e^{2i\t}, \, e^{-2i\t},
q e^{i\t}/d,  q e^{-i\t}/d; q\r)_\infty}
{\l(\a e^{i\t}, \a e^{-i\t}, q e^{i\t}/\a, q e^{-i\t}/\a; q\r)_\infty}
\, d \t, \notag
\end{align}
on the right side of (4.21). It is more convinient to write the $u_\nu$'s
instead of $v_\nu$'s here, cf. (2.14), (4.21), and (6.1).
Using (2.8) with $a$ and $c$ interchanged we get
\begin{align}
d^2_\nu
&=\frac{\l(abq^\nu, q^{1-\nu}/cd ; q\r)_\infty}
{(ab, q/cd; q)_\infty} \\
&\times\sum_{n=0}^\infty q^n\,
\frac{\l(q^{-\nu}, abcdq^{\nu-1}; q\r)_n}
{(q, ac, bc, cd; q)_n} \notag\\
&\times\int_0^{\pi}
\frac{u_{\nu}(\cos \t ; a, b, c; d)}
{\l(a e^{i\t}, a e^{-i\t}, b e^{i\t}, b e^{-i\t}, cq^n e^{i\t}, cq^n e^{-i\t}
; q\r)_\infty} \notag\\
&\qquad\times
\frac{\l(e^{2i\t}, \, e^{-2i\t},
q e^{i\t}/d,  q e^{-i\t}/d; q\r)_\infty}
{\l(\a e^{i\t}, \a e^{-i\t}, q e^{i\t}/\a, q e^{-i\t}/\a; q\r)_\infty}
\, d \t \notag\\
&+\frac{\l(q^{-\nu}, abcdq^{\nu-1}, qa/d, qb/d ; q\r)_\infty}
{(ab, ac, bc, cd/q; q)_\infty} \notag\\
&\times\sum_{n=0}^\infty q^n\,
\frac{\l(abq^\nu, q^{1-\nu}/cd; q\r)_n}
{(q, qa/d, qb/d, q^2/cd; q)_n} \notag\\
&\times\int_0^{\pi}
\frac{u_{\nu}(\cos \t ; a, b, c; d)}
{\l(a e^{i\t}, a e^{-i\t}, b e^{i\t}, b e^{-i\t},
q^{1+n} e^{i\t}/d,  q^{1+n} e^{-i\t}/d; q\r)_\infty} \notag\\
&\qquad\times
\frac{\l(e^{2i\t}, \, e^{-2i\t},
q e^{i\t}/d,  q e^{-i\t}/d; q\r)_\infty}
{\l(\a e^{i\t}, \a e^{-i\t}, q e^{i\t}/\a, q e^{-i\t}/\a; q\r)_\infty}
\, d \t. \notag
\end{align}
Both integrals on the right side have the same structure and can be
evaluated in a similar way. Let
\begin{align}
I_n(\g):=
&\int_0^{\pi}
\frac{u_{\nu}(\cos \t ; a, b, c; d)}
{\l(a e^{i\t}, a e^{-i\t}, b e^{i\t}, b e^{-i\t},
\g q^n e^{i\t}, \g q^n e^{-i\t}; q\r)_\infty} \\
&\qquad\times
\frac{\l(e^{2i\t}, \, e^{-2i\t},
q e^{i\t}/d,  q e^{-i\t}/d; q\r)_\infty}
{\l(\a e^{i\t}, \a e^{-i\t}, q e^{i\t}/\a, q e^{-i\t}/\a; q\r)_\infty}
\, d \t \notag\\
=&\frac{\l(qa/d, bcq^\nu ; q\r)_\infty}
{\l(q^{1-\nu}a/d, bc ; q\r)_\infty} \notag\\
&\times\int_0^{\pi}
\ _8W_7\l(\frac{a}{d}\,q^{-\nu};
q^{-\nu}, q^{1-\nu}/bd, q^{1-\nu}/cd,
ae^{i\t}, ae^{-i\t}; q,\, bcq^\nu \r) \notag\\
&\quad\times \frac{\l(e^{2i\t},\; e^{-2i\t}; q\r)_\infty}
{\l(a e^{i\t}, a e^{-i\t}, b e^{i\t}, b e^{-i\t},
\g q^n e^{i\t}, \g q^n e^{-i\t}; q\r)_\infty} \notag\\
&\qquad\times
\frac{\l(q^{1-\nu}e^{i\t}/d, \; q^{1-\nu}e^{-i\t}/d; q\r)_\infty}
{\l(\a e^{i\t}, \a e^{-i\t}, q e^{i\t}/\a, q e^{-i\t}/\a; q\r)_\infty}
\, d \t \notag
\end{align}
by (2.7) with $\g=c$ and $\g=q/d$. The last integral can be evaluated
in terms of two balanced $_5\f_4$'s by (6.4) of \cite{Is:Ra:Su},
\begin{align}
&\int_0^\pi
{}_8W_7\left(dg/q; h, r, g/f, de^{i\t}, de^{-i\t}; q, gf/hr\right) \\
&\quad\times
\frac{\l(e^{2i\t}, e^{-2i\t}; q\r)_\infty}
{\l(ae^{i\t}, ae^{-i\t}, be^{i\t}, be^{-i\t},
ce^{i\t}, ce^{-i\t}, de^{i\t}, de^{-i\t}; q\r)_\infty} \notag\\
&\qquad\qquad\qquad\qquad\times
\frac{\l(ge^{i\t}, ge^{-i\t}; q\r)_\infty}
{\l(fe^{i\t}, fe^{-i\t}; q\r)_\infty}
\, d\t \notag\\
&= \frac{2\pi \,(abcd, dg, g/d;q)_\infty}
{(q, ab, ac, ad, bc, bd, cd, df, f/d;q)_\infty}  \notag\\
&\quad\times
\ _5\phi_4\left(
\begin{array}{cc}
ad, bd, cd, g/f, dg/hr \\
\, \\
abcd, qd/f, dg/h, dg/r
\end{array}
; q,\, q \right) \notag\\
&+\frac{2\pi \,(abcf, dg, g/f;q)_\infty}
{(q, ab, ac, af, bc, bf, cf, df, d/f ;q)_\infty} \notag\\
&\quad\times
\frac{(fg/h, fg/r, dg/hr;q)_\infty}{(dg/h, dg/r, fg/hr;q)_\infty}
\notag\\
&\quad \times
\ _5\phi_4\left(
\begin{array}{cc}
af, bf, cf, g/d, gf/hr \\
\, \\
abcf, qf/d, gf/h, gf/r
\end{array}
; q, \,q \right), \notag
\end{align}
this formula appears in a straightforward manner if one uses
Bailey's transform (III.36) of \cite{Ga:Ra} and
the Askey--Wilson integral. Thus,
\begin{align}
I_n(\g)=
&\frac{2\pi\l(qa/d, q^{1-\nu}/ad, bcq^\nu, a\g q;q\r)_\infty}
{\l(q, q, ab, bc, b/a, \a a, qa/\a, \a \g, q\g/\a, a\g; q\r)_\infty}\\
&\times\frac{\l(\a \g, q\g/\a, a\g ;q\r)_n}{\l(a\g q;q\r)_n}\notag\\
&\times
\ _5\phi_4\left(
\begin{array}{cc}
a\g q^n, acq^\nu, q^{1-\nu}/bd, \a a, qa/\a \\
\, \\
a\g q^{1+n}, ac, qa/b, qa/d
\end{array}
; q,\, q \right) \notag\\
&+\frac{2\pi\l(qb/d, q^{1-\nu}/bd, acq^\nu, b\g q;q\r)_\infty}
{\l(q, q, ab, ac, a/b, \a b, qb/\a, \a \g, q\g/\a, b\g; q\r)_\infty}\notag\\
&\times\frac{\l(\a \g, q\g/\a, b\g ;q\r)_n}{\l(b\g q;q\r)_n}\notag\\
&\times
\ _5\phi_4\left(
\begin{array}{cc}
b\g q^n, bcq^\nu, q^{1-\nu}/ad, \a b, qb/\a \\
\, \\
b\g q^{1+n}, bc, qb/a, qb/d
\end{array}
; q,\, q \right), \notag
\end{align}
where $\g=c$ and $\g=q/d$. One can see that the second term here
equals the first one with $a$ and $b$ interchanged.

Combining (6.2), (6.3), and (6.5), we obtain
\begin{align}
d_\nu^2
&=\int_0^{\pi}
\frac{v_{\nu}^2(\cos \t ; a, b, c; d)}
{\l(a e^{i\t},  a e^{-i\t}, b e^{i\t}, b e^{-i\t}, c e^{i\t},  c e^{-i\t},
q e^{i\t}/d,  q e^{-i\t}/d; q\r)_\infty}\\
&\qquad\times
\frac{\l(e^{2i\t}, \, e^{-2i\t}; q\r)_\infty}
{\l(\a e^{i\t}, \a e^{-i\t}, q e^{i\t}/\a, q e^{-i\t}/\a; q\r)_\infty}
\, d \t \notag \\
&=\frac{2\pi}{(1-ac)(q, ab; q)^2}\,
\frac
{\l(qa/d, abq^\nu, bcq^\nu, q^{1-\nu}/ad, q^{1-\nu}/cd;q\r)_\infty}
{\l(bc, b/a, q/cd, \a a, qa/\a, \a c, qc/\a;q\r)_\infty} \notag\\
&\quad\times\sum_{n=0}^\infty q^n\,
\frac{\l(q^{-\nu}, abcdq^{\nu-1}, \a c, qc/\a; q\r)_n}
{(q, \a cq, bc, cd; q)_n} \notag\\
&\quad\times
\ _5\phi_4\left(
\begin{array}{cc}
acq^n, acq^\nu, q^{1-\nu}/bd, \a a, qa/\a \\
\, \\
acq^{n+1}, ac, qa/b, qa/d
\end{array}
; q,\, q \right) \notag\\
&+\frac{2\pi}{(1-bc)(q, ab; q)^2}\,
\frac
{\l(qb/d, abq^\nu, acq^\nu, q^{1-\nu}/bd, q^{1-\nu}/cd;q\r)_\infty}
{\l(ac, a/b, q/cd, \a b, qb/\a, \a c, qc/\a;q\r)_\infty} \notag\\
&\quad\times\sum_{n=0}^\infty q^n\,
\frac{\l(q^{-\nu}, abcdq^{\nu-1}, \a c, qc/\a; q\r)_n}
{(q, q\a c, ac, cd; q)_n} \notag\\
&\quad\times
\ _5\phi_4\left(
\begin{array}{cc}
bcq^n, bcq^\nu, q^{1-\nu}/ad, \a b, qb/\a \\
\, \\
bcq^{n+1}, bc, qb/a, qb/d
\end{array}
; q,\, q \right) \notag\\
&+\frac{2\pi\l(qa/d;q\r)_\infty^2}{(1-qa/d)(q, ab, bc; q)^2}\,
\frac
{\l(q^{-\nu}, abcdq^{\nu-1}, q^{1-\nu}/ad, bcq^\nu, qb/d;q\r)_\infty}
{\l(ac, b/a, cd/q, \a a, qa/\a, q\a/d, q^2/\a d;q\r)_\infty} \notag\\
&\quad\times\sum_{n=0}^\infty q^n\,
\frac{\l(abq^\nu, q^{1-\nu}/cd, q\a/d, q^2/\a d; q\r)_n}
{(q, qb/d, q^2a/d, q^2/cd; q)_n} \notag\\
&\quad\times
\ _5\phi_4\left(
\begin{array}{cc}
aq^{n+1}/d, acq^\nu, q^{1-\nu}/bd, \a a, qa/\a \\
\, \\
aq^{n+2}/d, ac, qa/b, qa/d
\end{array}
; q,\, q \right) \notag\\
&+\frac{2\pi\l(qb/d;q\r)_\infty^2}{(1-qb/d)(q, ab, ac; q)^2}\,
\frac
{\l(q^{-\nu}, abcdq^{\nu-1}, q^{1-\nu}/bd, acq^\nu, qa/d;q\r)_\infty}
{\l(bc, a/b, cd/q, \a b, qb/\a, q\a/d, q^2/\a d;q\r)_\infty} \notag\\
&\quad\times\sum_{n=0}^\infty q^n\,
\frac{\l(abq^\nu, q^{1-\nu}/cd, q\a/d, q^2/\a d; q\r)_n}
{(q, qa/d, q^2b/d, q^2/cd; q)_n} \notag\\
&\quad\times
\ _5\phi_4\left(
\begin{array}{cc}
bq^{n+1}/d, bcq^\nu, q^{1-\nu}/ad, \a b, qb/\a \\
\, \\
bq^{n+2}/d, bc, qb/a, qb/d
\end{array}
; q,\, q \right). \notag
\end{align}
One can see again that the second and the forth terms in this formula
are equal to the first and the third ones, respectively,
with $a$ and $b$ interchanged.
When $\nu$ satisfies the boundary condition (4.22) the last integral
gives the values of the normalization constants in the orthogonality
relation (4.23)--(4.24).

\section{Some Identity}

In this section we shall derive an interesting relation
involving a determinant of four $_8\f_7$-functions of type (2.7).
Let
\begin{align}
&u(z)=u_\nu(x(z); a, b, c; d),\\
&v(z)=u_\nu(x(z); a, b, d; c)=\l . u(z)\r|_{c\leftrightarrow d}
\end{align}
be two solutions of equation (2.1) corresponding to
the same eigenvalue (2.4). Then, due to (4.3),
\begin{equation}
\Delta \l[\sigma(z)\, \rho(z)\; W\l(u(z),\; v(z)\r)\r]=0,
\end{equation}
where $W(u, v)$ is the analog of the Wronskian defined in (4.4) and
$\rho(z)$ is the appropriate solution of the Pearson equation (3.7).
Using the difference-differentiation formula (2.15) we can rewrite
the ``Wronskian" as
\begin{align}
W\l(u(z),\; v(z)\r)=
&\frac{2q}{(1-q)c}\,
\frac{\l(1-q^{-\nu}\r)\l(1-abcdq^{\nu-1}\r)}
{(1-ab)(1-ad)(1-bd)} \\
&\times u_{\nu-1}\l(x(z-1/2); aq^{1/2}, bq^{1/2}, dq^{1/2}; cq^{1/2}\r)
\notag\\
&\times u_\nu\l(x(z); a, b, c; d\r)\notag\\
&-\frac{2q}{(1-q)d}\,
\frac{\l(1-q^{-\nu}\r)\l(1-abcdq^{\nu-1}\r)}
{(1-ab)(1-ac)(1-bc)} \notag \\
&\quad\times u_{\nu-1}\l(x(z-1/2); aq^{1/2}, bq^{1/2}, cq^{1/2}; dq^{1/2}\r)
\notag\\
&\quad\times u_\nu\l(x(z); a, b, d; c\r).
\notag
\end{align}
One can easily see that the function
\begin{align}
g(z)
&=\sigma(z)\, \rho(z)\; W\l(u(z),\; v(z)\r)\\
&=cd\,\frac{\l(q^z/c, q^{1-z}/c, q^z/d, q^{1-z}/d ;q\r)_\infty}
{\l(aq^z, aq^{1-z}, bq^z, bq^{1-z};q\r)_\infty} \notag\\
&\qquad\qquad\times W\l(u(z),\; v(z)\r) \notag
\end{align}
in (7.3) is doubly periodic function in $z$ without poles
in the rectangle on the Figure. Therefore, this function is just
a constant by Liouville's theorem,
\begin{equation}
\sigma(z)\, \rho(z)\; W\l(u(z),\; v(z)\r)=C.
\end{equation}
To find the value of this constant we can choose here $z=z_0$ such that
$q^{z_0}=a$.
From (7.4) and (2.7) one gets
\begin{align}
&W\l(u(z_0),\; v(z_0)\r)\\
&=\frac{2q\l(1-q^{-\nu}\r)\l(1-abcdq^{\nu-1}\r)}
{(1-q)} \notag \\
&\times\Bigg[
\frac{\l(q, qbc, qb/d, qc/d, bdq^\nu,
q^{1-\nu}/ac, q^{1-\nu}a/d, a^2bcq^\nu ; q\r)_\infty}
{acd\l(ab, qac, bc, bd, 1/ac, q/ad, a/d, qabc/d; q\r)_\infty} \notag\\
&\qquad\times
\ _8W_7\l(abc/d; ab, ac, a/d, q^{1-\nu}/ad, bcq^\nu ; q, \; q\r) \notag\\
&\quad-\frac{\l(q, qbd, qb/c, qd/c, bcq^\nu,
q^{1-\nu}/ad, q^{1-\nu}a/c, a^2bdq^\nu ; q\r)_\infty}
{acd\l(ab, qad, bd, bc, 1/ad, q/ac, a/c, qabd/c; q\r)_\infty} \notag\\
&\quad\qquad\times
\ _8W_7\l(abd/c; ab, ad, a/c, q^{1-\nu}/ac, bdq^\nu ; q, \; q\r)
\Bigg] \notag\\
&\quad=\frac
{2q\l(q, a^2, qb/a, d/c, qc/d,
q^{-\nu}, abcdq^{\nu-1},  abq^\nu, q^{1-\nu}/cd; q\r)_\infty}
{(1-q)acd\l(ab, qac, ad, bc, bd, 1/ac, q/ad, a/c, a/d; q\r)_\infty} \notag
\end{align}
by (III.24) and (II.25) of \cite{Ga:Ra}.

As a result, from (7.6) and (7.7) we find the value of the ``Wronskian" of
the $_8\f_7$-functions (7.1) and (7.2),
\begin{align}
W\l(u(z),\; v(z)\r)=
&\frac{2q\l(c/d, qd/c, q^{-\nu}, abcdq^{\nu-1},
abq^\nu, q^{1-\nu}/cd; q\r)_\infty}
{(1-q)c\l(ab, ab, ac, ad, bc, bd; q\r)_\infty} \\
&\times\frac{\l(aq^z, aq^{1-z}, bq^z, bq^{1-z};q\r)_\infty}
{\l(q^z/c, q^{1-z}/c, q^z/d, q^{1-z}/d ;q\r)_\infty}.
\notag
\end{align}
Due to (7.4) and (2.7) the last equation can also be rewritten
in a more explicit form,
\begin{align}
&\frac{d \l(ac, adq; q\r)_\infty}
{\l(q^{2-\nu}a/c, q^{1-\nu}a/d; q\r)_\infty}\\
&\times
\frac{\l(q^{1-\nu+z}/c, q^{2-\nu-z}/c,
q^{1-\nu+z}/d, q^{1-\nu-z}/d; q\r)_\infty}
{\l(q^z/c, q^{1-z}/c, q^{1+z}/d, q^{1-z}/d; q\r)_\infty}\notag\\
&\times
\ _8W_7\l(q^{1-\nu}a/c; q^{1-\nu}, q^{1-\nu}/bc, q^{1-\nu}/cd,
aq^z, aq^{1-z}; q, \; q\r) \notag\\
&\times
\ _8W_7\l(q^{-\nu}a/d; q^{-\nu}, q^{1-\nu}/bd, q^{1-\nu}/cd,
aq^z, aq^{-z}; q, \; q\r) \notag\\
&-\frac{c \l(ad, acq; q\r)_\infty}
{\l(q^{2-\nu}a/d, q^{1-\nu}a/c; q\r)_\infty} \notag\\
&\quad\times
\frac{\l(q^{1-\nu+z}/d, q^{2-\nu-z}/d,
q^{1-\nu+z}/c, q^{1-\nu-z}/c; q\r)_\infty}
{\l(q^z/d, q^{1-z}/d, q^{1+z}/c, q^{1-z}/c; q\r)_\infty}\notag\\
&\quad\times
\ _8W_7\l(q^{1-\nu}a/d; q^{1-\nu}, q^{1-\nu}/bd, q^{1-\nu}/cd,
aq^z, aq^{1-z}; q, \; q\r) \notag\\
&\quad\times
\ _8W_7\l(q^{-\nu}a/c; q^{-\nu}, q^{1-\nu}/bc, q^{1-\nu}/cd,
aq^z, aq^{-z}; q, \; q\r) \notag\\
&=\frac{\l(c/d, qd/c, q^{1-\nu}, abcdq^\nu,
abq^\nu, q^{1-\nu}/cd; q\r)_\infty}
{c\l(ab, abq, qa/c, qa/d, bcq^\nu, bdq^\nu; q\r)_\infty}\notag\\
&\quad\times\frac{d \l(aq^z, aq^{1-z}, bq^z, bq^{1-z};q\r)_\infty}
{\l(q^z/c, q^{1-z}/c, q^z/d, q^{1-z}/d ;q\r)_\infty}.
\notag
\end{align}
The second term on the left side here is the same as the first one
with $c$ and $d$ interchanged.

One can see from (7.8) that two solutions
$u(z)$ and $v(z)$
are linear dependent when $\nu$ is an integer.
In this case due to (2.13) both solutions are
the Askey--Wilson polynomials, up to a factor,
which are related by Sears' transformation.
On the other hand, equation (7.8) shows that there is no analog of
Sears's transformation at the level of
very-well-poised $_8\phi_7$-functions.

\section{Some Special and Limiting Cases}

The Askey--Wilson polynomials are known as the most general system of
classical orthogonal polynomials. They include all other classical
orthogonal polynomials as special and/or limiting cases
\cite{An:As}, \cite{As:Wi}, \cite{Ga:Ra}, \cite{Ko:Sw}, and \cite{Ni:Su:Uv}.
Let us discuss in a similar manner a few interesting special cases of
the orthogonal $_8\f_7$-functions (2.7).

\subsection{Extension of continuous dual $q$-Hahn polynomials\/}

Letting $c\to 0$ in (2.7) and then changing $d$ by $c$ we get
\begin{align}
u_{\nu}(x; a, b; c)
&=\frac
{\l(qa/c, q^{1-\nu+z}/c, q^{1-\nu-z}/c; q\r)_\infty}
{\l(q^{1-\nu}a/c, q^{1+z}/c, q^{1-z}/c; q\r)_\infty} \notag\\
&\quad\times \ _7\f_7\l(
\begin{array}{cc}
\dfrac{aq^{-\nu}}{c},
q \sqrt{\dfrac{aq^{-\nu}}{c}},-q \sqrt{\dfrac{aq^{-\nu}}{c}},
q^{-\nu}, \dfrac{q^{1-\nu}}{bc}, a q^z, aq^{-z} \\
\, \\
\sqrt{\dfrac{aq^{-\nu}}{c}}, -\sqrt{\dfrac{aq^{-\nu}}{c}},
\dfrac{aq}{c}, ab, 0, \dfrac{q^{1-\nu+z}}{c}, \dfrac{q^{1-\nu-z}}{c}
\end{array}
;\; q, \; q\,\frac{b}{c} \r) \notag\\
&=\frac{\l(q^{1-\nu}/ac;q \r)_\infty}
{\l(q/ac;q \r)_\infty}\;
\ _3\f_2\l(
\begin{array}{cc}
q^{-\nu}, a q^z, aq^{-z} \\
\, \\
ab, ac
\end{array}
;\; q, \; q \r) \notag\\
&\quad+
\frac{\l(q^{-\nu}, qb/c, aq^z, aq^{-z}; q\r)_\infty}
{\l(ab, ac/q, q^{1+z}/c, q^{1-z}/c; q\r)_\infty} \notag\\
%\;
&\qquad\times
\ _3\f_2\l(
\begin{array}{cc}
q^{1-\nu}/ac, q^{1+z}/c, q^{1-z}/c \\
\, \\
qb/c, q^2/ac
\end{array}
;\; q, \; q \r).
\end{align}
The $_7\f_7$-function here can also be transformed to a $_3\f_2$
by (3.2.11) of \cite{Ga:Ra},
\begin{align}
u_{\nu}(x; a, b; c)
=&\frac{\l(abq^\nu, q^{1-\nu+z}/c, q^{1-\nu-z}/c; q \r)_\infty}
{\l(ab, q^{1+z}/c, q^{1-z}/c ; q \r)_\infty}\\
&\times \ _3\f_2\l(
\begin{array}{cc}
q^{-\nu},  q^{1-\nu}/ac, q^{1-\nu}/bc \\
\, \\
q^{1-\nu+z}/c, q^{1-\nu-z}/c
\end{array}
;\; q, \; abq^\nu \r).
\notag
\end{align}
One can easily see that for an integer $\nu$ our function (8.1)--(8.2)
is just a multiple of the continuous dual $q$-Hahn polynomial
\cite{An:As}, \cite{As:Wi}, \cite{Ko:Sw}, and \cite{Ni:Su:Uv}.

The orthogonality relation for the corresponding entire function,
\begin{equation}
v_\nu(x; a, b; c)=\l(q^{1+z}/c, q^{1-z}/c; q\r)_\infty\;
u_\nu(x; a, b; c),
\end{equation}
takes the form
\begin{align}
&\int_0^{\pi}
\frac{v_{\mu}(\cos \t ; a, b; c) \;
v_{\nu}(\cos \t ; a, b; c)}
{\l(a e^{i\t},  a e^{-i\t}, b e^{i\t}, b e^{-i\t},
q e^{i\t}/c,  q e^{-i\t}/c; q\r)_\infty}\\
&\qquad\times
\frac{\l(e^{2i\t}, \, e^{-2i\t}; q\r)_\infty}
{\l(\a e^{i\t}, \a e^{-i\t}, q e^{i\t}/\a, q e^{-i\t}/\a; q\r)_\infty}
\, d \t =0 \notag
\end{align}
if $\mu \ne \nu$, and
\begin{align}
&\int_0^{\pi}
\frac{\l(v_{\nu}(\cos \t ; a, b; c)\r)^2}
{\l(a e^{i\t},  a e^{-i\t}, b e^{i\t}, b e^{-i\t},
q e^{i\t}/c,  q e^{-i\t}/c; q\r)_\infty}\\
&\qquad\times
\frac{\l(e^{2i\t}, \, e^{-2i\t}; q\r)_\infty}
{\l(\a e^{i\t}, \a e^{-i\t}, q e^{i\t}/\a, q e^{-i\t}/\a; q\r)_\infty}
\, d \t \notag \\
&\quad=\l(q, q, \a a, qa/\a, \a b, qb/\a ,
q\a/c, q^2/{\a c}; q\r)_\infty^{-1} \notag \\
&\qquad\times
\frac{-4\pi q^{1/2} c}{1-q} \;\;
\frac{\partial}{\partial \lambda_\nu}\,
v_\nu(\eta; a, b; c)\;
\l(\frac{\nabla}{\nabla x}\,
v_\nu(x; a, b; c)\r)\Bigg|_{x=\eta} \notag
\end{align}
if $\mu=\nu$, respectively.
The ``degree" parameters $\mu$ and $\nu$ here satisfy
the ``boundary" condition
\begin{align}
&v_{\ep}\l(\dfrac{1}{2}\l(\a+\a^{-1}\r); a, b; c\r)\\
&=\l(q\a/c, q^{1-\varepsilon}/\a c; q\r)_\infty
\ _3\f_2\l(
\begin{array}{cc}
q^{-\varepsilon},  a/\a, b/\a\\
\, \\
q^{1-\varepsilon}/\a c, ab
\end{array}
;\; q, \; q\,\frac{\a}{c}\r)
=0.\notag
\end{align}
Properties of zeros of this function follow as a special case
of the results of Section 5.

\subsection{Extension of Al-Salam and Chihara polynomials\/}

Letting $b\to 0$ in (8.1)--(8.2) and then changing $c$ by $b$ one gets
\begin{align}
u_{\nu}(x; a; b)
=&\frac{\l(q^{1-\nu+z}/b, q^{1-\nu-z}/b; q \r)_\infty}
{\l(q^{1+z}/b, q^{1-z}/b ; q \r)_\infty}\\
&\times \ _2\f_2\l(
\begin{array}{cc}
q^{-\nu},  q^{1-\nu}/ab \\
\, \\
q^{1+z}/b, q^{1-z}/b
\end{array}
;\; q, \; q\,\frac{a}{b}\r) \notag\\
=&\frac{\l(q^{1-\nu}/ab;q \r)_\infty}
{\l(q/ab;q \r)_\infty}\;
\ _3\f_2\l(
\begin{array}{cc}
q^{-\nu}, a q^z, aq^{-z} \\
\, \\
ab, 0
\end{array}
;\; q, \; q \r) \notag\\
&+
\frac{\l(q^{-\nu}, aq^z, aq^{-z}; q\r)_\infty}
{\l(ab/q, q^{1+z}/b, q^{1-z}/b; q\r)_\infty} \notag\\
&\quad\times
\ _3\f_2\l(
\begin{array}{cc}
q^{1-\nu}/ab, q^{1+z}/b, q^{1-z}/b \\
\, \\
q^2/ab, 0
\end{array}
;\; q, \; q \r). \notag
\end{align}
For an integer $\nu$ this function is a multiple of
the Al-Salam and Chihara polynomial
\cite{An:As}, \cite{As:Ra:Su}, \cite{As:Wi}, and \cite{Ko:Sw}.

The orthogonality relation of the entire function,
\begin{equation}
v_\nu(x; a; b)=\l(q^{1+z}/b, q^{1-z}/b; q\r)_\infty\;
u_\nu(x; a; b),
\end{equation}
is
\begin{align}
&\int_0^{\pi}
\frac{v_{\mu}(\cos \t ; a; b) \;
v_{\nu}(\cos \t ; a; b)}
{\l(a e^{i\t},  a e^{-i\t},
q e^{i\t}/b,  q e^{-i\t}/b; q\r)_\infty}\\
&\times
\frac{\l(e^{2i\t}, \, e^{-2i\t}; q\r)_\infty}
{\l(\a e^{i\t}, \a e^{-i\t}, q e^{i\t}/\a, q e^{-i\t}/\a; q\r)_\infty}
\, d \t =0, \notag
\end{align}
where $\mu \ne \nu$ are solutions of
\begin{align}
&v_{\ep}\l(\dfrac{1}{2}\l(\a+\a^{-1}\r); a; b\r)\\
&=\l(q\a/b, q^{1-\varepsilon}/\a b; q\r)_\infty
\ _2\f_1\l(
\begin{array}{cc}
q^{-\varepsilon},  a/\a\\
\, \\
q^{1-\varepsilon}/\a b
\end{array}
;\; q, \; q\,\frac{\a}{b}\r)
=0.\notag
\end{align}
For properties of zeros see Section 5.

\subsection{Extension of continuous big $q$-Hermite polynomials\/}

Letting $a\to 0$ in (8.7) and then changing $b$ by $a$ we have
\begin{align}
u_{\nu}(x; a)
=&\frac{\l(q^{1-\nu+z}/a, q^{1-\nu-z}/a; q \r)_\infty}
{\l(q^{1+z}/a, q^{1-z}/a ; q \r)_\infty} \\
&\times \ _1\f_2\l(
\begin{array}{cc}
q^{-\nu} \\
\, \\
q^{1-\nu+z}/a, q^{1-\nu-z}/a
\end{array}
;\; q, \; q^{2-\nu}/a^2 \r).
\notag
\end{align}
For an integer $\nu$ this function is a multiple of
the continuous big $q$-Hermite polynomial \cite{Ko:Sw}.

The orthogonality relation has the form
\begin{equation}
\int_0^{\pi}
\frac{v_{\mu}(\cos \t ; a) \; v_{\nu}(\cos \t ; a)\;
\l(e^{2i\t}, \, e^{-2i\t}; q\r)_\infty}
{\l(
\a e^{i\t}, \a e^{-i\t}, q e^{i\t}/\a, q e^{-i\t}/\a,
q e^{i\t}/a,  q e^{-i\t}/a ; q\r)_\infty}
\, d \t =0,
\end{equation}
where $\mu \ne \nu$ satisfy
\begin{align}
&v_{\ep}\l(\dfrac{1}{2}\l(\a+\a^{-1}\r); a\r)\\
&=\l(q^{1-\varepsilon}/\a a; q\r)_\infty
\ _1\f_1\l(
\begin{array}{cc}
q/\a a\\
\, \\
q^{1-\varepsilon}/\a a
\end{array}
;\; q, \; q^{1-\varepsilon}\a/a\r)
=0 \notag
\end{align}
and $v_{\varepsilon}(x; a)=\l(q^{1+z}/a, q^{1-z}/a; q\r)_\infty\;
u_{\varepsilon}(x; a)$.

\subsection{Extension of continuous $q$-Hermite polynomials}

The continuous $q$-Hermite polynomials $H_n(x|q)$
are the simplest special case $a=b=c=d=0$ of
the Askey--Wilson polynomials $p_n(x; a, b, c, d)$
or the special case $a=0$ of
the continuous big $q$-Hermite polynomials $p_n(x; a)$
\cite{As:Wi}, \cite{Ga:Ra}, \cite{Ko:Sw}, and \cite{Ni:Su:Uv}.
See \cite{As0} and \cite{As:Ra:Su} for the proof of
$$\lim_{a \to 0} p_n(x; a)=H_n(x|q)$$
directly from the series representation of
the continuous big $q$-Hermite polynomials.

Let us consider the difference equation (2.1) with $a=b=c=d=0$ and
let us choose the following solution,
$u_\nu(z)=H_\nu\l(x(z)|q\r)$, such that
\begin{align}
H_{\nu}(x|q)
=&\frac{\l(-q^{1-\nu+2z}, -q^{1-\nu-2z}; q^2 \r)_\infty}
{\l(-q^{1+2z}, -q^{1-2z}; q^2 \r)_\infty} \\
&\times \ _2\f_2\l(
\begin{array}{cc}
q^{-\nu}, \quad q^{1-\nu} \\
\, \\
-q^{1-\nu+2z},\; -q^{1-\nu-2z}
\end{array}
;\; q^2, \; q\r)
\notag
\end{align}
as a nonterminating extension of the continuous $q$-Hermite polynomials.
This solution differs from the corresponding one in \cite{At:Su}
by a periodic factor. For an integer $\nu$
function (8.14) coinsides with $H_n(x|q)$ up to a constant.

Comparing (8.14) with the equation (2.8) of \cite{Bu:Su}
one can see that function $H_\nu(x|q)$ is a multiple of
the basic cosine function $C(x; \omega)$ for certain values of
parameter $\nu$. Therefore, this function satisfies the orthogonality
relation (1.16) under the ``boundary" conditions (1.20).
We would like to leave the details to the reader.

In a similar fashion one can consider some other special and limiting cases
of the orthogonal $_8\f_7$-functions.

\section{Appendix: Estimate of Number of Zeros}

In this section we give an estimate for number of zeros of the ``boundary"
function $v_\nu(\eta ;a, b, c; d)$ on the basis of Jensen's theorem
(see, for example, \cite{Lev}). We shall apply the method proposed by
Mourad Ismail at the level of the third Jackson $q$-Bessel functions
\cite{Is} (see also \cite{Is:Ma:Su3} and \cite{Bu:Su}
for an extension of his idea to $q$-Bessel functions on a $q$-quadratic grid
and $q$-trigonometric functions, respectively).

Let us consider the entire function
\begin{align}
f(\zeta)
&=v_\nu\l(\frac{1}{2}\l(\a + \a^{-1}\r); a, b, c; d\r) \\
&=\frac{\l(q\a/d ; q\r)_\infty}
{\l(abc/\a ; q\r)_\infty} \notag \\
&\quad\times\prod_{k=0}^\infty
\l(1-\zeta q^{k+1}/\a d + abcq^{2k+1}/\a^2 d \r) \notag \\
&\quad\times
\sum_{m=0}^\infty
\l(q\frac{\a}{d}\r)^m
\frac{\l(1-abcq^{2m-1}/\a \r) \l(abcq^{-1}/\a, a/\a , b/\a , c/\a ; q\r)_m}
{\l(1-abcq^{-1}/\a \r) \l(q, ab, ac, bc; q\r)_m} \notag \\
&\qquad\times\prod_{k=0}^{m-1}\frac
{1-\zeta q^k + abcdq^{2k-1}}
{1-\zeta q^{k+1}/\a d + abcq^{2k+1}/\a^2 d}
\notag
\end{align}
in a complex variable
\begin{equation}
\zeta=q^{-\nu}+abcdq^{\nu-1},
\end{equation}
where $|q\a/d|<1$ (cf. (5.2)).

Let $n_{f}\left( r\right) $ be the number of of zeros of $f(\zeta )$ in the
circle $\left| \zeta \right| <r$. Consider also circles of radius
$R=R_{n}=q^{-n}/\b + \b abcdq^{n-1}$, $q< \b < 1$ with
$n=1,2,3, \ldots$ in the complex $\zeta $-plane.
Since $n_{f}\left( r\right) $ is nondecreasing with $r$ one can write
\begin{equation}
n_{f}\left( R_{n}\right) \leq n_{f}\left( r\right) \leq n_{f}\left(
R_{n+1}\right)
\end{equation}
if $R_{n}\leq r\leq R_{n+1}$, and, therefore,
\begin{equation}
n_{f}\left( R_{n}\right) \ \int_{R_{n}}^{R_{n+1}}\frac{dr}{r}\leq
\int_{R_{n}}^{R_{n+1}}\frac{n_{f}\left( r\right) }{r}\ dr\leq n_{f}\left(
R_{n+1}\right) \ \int_{R_{n}}^{R_{n+1}}\frac{dr}{r}.
\end{equation}
But
\begin{equation}
\left. \int_{R_{n}}^{R_{n+1}}\frac{dr}{r}=\log r\right|
_{R_{n}}^{R_{n+1}}=\log
\l(\frac{q^{-1}+\b^2 abcdq^{2n}}{1+\b^2 abcdq^{2n-1}}\r)
=\log q^{-1} + \text{o}(1)
\end{equation}
as $n\to\infty$,
and, finally, one gets
\begin{equation}
\log q^{-1}\ n_{f}\left( R_{n}\right) < \int_{R_{n}}^{R_{n+1}}
\frac
{n_{f}\left( r\right) }{r}\ dr < \log q^{-1}\ n_{f}\left( R_{n+1}\right)
\end{equation}
for sufficiently large $n$.

The next step is to estimate the integral in (9.6). By Jensen's
theorem \cite{Lev}
\begin{align}
\int_{R_{n}}^{R_{n+1}}\frac{n_{f}\left( r\right) }{r}\ dr
&=\int_{0}^{R_{n+1}}\frac{n_{f}\left( r\right) }{r}\ dr-\int_{0}^{R_{n}}%
\frac{n_{f}\left( r\right) }{r}\ dr \\
&=\frac{1}{2\pi }\int_{0}^{2\pi }\log \left|
\frac{f\left( R_{n+1}\, e^{i\vartheta }\right) }
{f\left( R_n\, e^{i\vartheta}\right) }\right| \ d\vartheta .  \notag
\end{align}
For large values of $n$ we have $R_n \sim q^{-n}/\b$ and,
in view of (5.3) and (9.1),
\begin{align}
\frac{f\l( R_{n+1}\, e^{i\vartheta}\r)}
{f\l( R_n\, e^{i\vartheta}\r)}
&\sim\prod_{k=0}^\infty\frac
{1-e^{i\vartheta}\,q^{k-n}/\a \b d + abcq^{2k+1}/\a^2 d}
{1-e^{\vartheta}\,q^{k-n+1}/\a \b d + abcq^{2k+1}/\a^2 d} \\
&\sim -e^{i\vartheta}q^{-n}/\a\b d,
\notag
\end{align}
where we have used the formula
\begin{align}
&\prod_{k=0}^\infty
\l(1-hq^{k-n}+gq^{2k}\r) \\
&=\prod_{k=0}^{n-1}
\l(1-hq^{k-n}+gq^{2k}\r)\;
\prod_{k=n}^\infty
\l(1-hq^{k-n}+gq^{2k}\r) \notag \\
&=\l(-h\r)^n q^{-n(n+1)/2}
\prod_{k=0}^{n-1}
\l(1-\frac{1}{h}\,q^{k+1}-\frac{g}{h}\,q^{2n-k-1}\r) \notag \\
&\qquad\qquad\qquad\qquad\times\prod_{k=0}^\infty
\l(1-hq^{k}+gq^{2n+2k}\r) \notag
\end{align}
with $h=qe^{i\vartheta}/\a \b d$ and $g=qabc/\a^2 d$.
Thus,
\begin{equation*}
\log \left| \frac{f\left( R_{n+1}\, e^{i\vartheta }\right) }
{f\left(R_n\, e^{i\vartheta }\right) }\right| \sim n\ \log q^{-1}
-\log \gamma,
\end{equation*}
where $\gamma = |\a \b d|$, and
\begin{equation}
\int_{R_{n}}^{R_{n+1}}\frac{n_{f}\left( r\right) }{r}\, dr
=n\log q^{-1}-\log \gamma +\text{o}\left( 1\right)
\end{equation}
as $n\rightarrow \infty $.

From (9.6) and (9.10),
\begin{equation}
1-\frac{\log \gamma /\log q^{-1}}{n}-\frac{1}{n}
< \frac{n_{f}\l(R_{n}\r)}{n}
< 1-\frac{\log \gamma /\log q^{-1}}{n}
\end{equation}
and, therefore,
\begin{equation}
\lim_{n\rightarrow \infty }\frac{n_{f}\left( R_{n}\right) }{n}=1.
\end{equation}
On the other hand, from (9.11),
\begin{equation}
n-1-\log \gamma /\log q^{-1}
< n_{f}\l(R_{n}\r)
< n-\log \gamma /\log q^{-1}.
\end{equation}
The difference between the upper and the lower bounds here is $1$
which means that there is only one positive root of
$v_\nu(\eta ; a, b, c; d)$ between the test points
$\nu=\omega_n$ and $\nu=\omega_{n+1}$ defined on the page~15
during the proof of Theorem 5.1 for large values of $n$.

\section*{Acknowledgments}

I would like to thank Dick Askey for encouraging conversations.
I wish to thank Joaqu{\'{\i}}n Bustoz, George Gasper,
Alberto Gr\"{u}nbaum,
Mourad Ismail, and Mizan Rahman for valuable discussions and comments.
This paper was prepaired while the author visited
Mathematical Sciences Research Institute in Berkeley, California,
a beautiful place to work,
and I gratefully acknowledge their hospitality.

\end{document}